\documentclass[final]{amsart}
\usepackage{pdfsync}
\usepackage{amsmath,amssymb}
\usepackage{enumerate}
\usepackage{xspace}
\usepackage{boxedminipage}
\usepackage{algorithmicx}
\usepackage[ruled]{algorithm}
\usepackage{algpseudocode}
\usepackage{listings}
\usepackage{tabularx}
\usepackage{subfigure}
\usepackage{tristan}
\usepackage{fullpage}

\usepackage[backend=bibtex,style=numeric,doi=false,eprint=false,url=false,isbn=false,bibstyle=numeric,sorting=nty,minbibnames=5,maxbibnames=5,defernumbers=true,firstinits=true]{biblatex}

\renewcommand{\bi}[2]{\ensuremath{\cA\qp{#1,#2}}}

\renewcommand{\vec}[1]{\geovec{#1}}
\renewcommand{\mat}[1]{\geomat{#1}}

\renewcommand{\enorm}[2]{\ensuremath{\Norm{#1}}_{dG,{#2}}}

\renewcommand{\eenorm}[2]{\,|\!|\!| {{#1}} |\!|\!|_{dG,{#2}}}
\newcommand{\sobg}[2]{\ensuremath{\WW^{#1,#2}_g}}
\renewcommand{\sobz}[2]{\ensuremath{\WW^{#1,#2}_0}}
\renewcommand{\sob}[2]{\ensuremath{\WW^{#1,#2}}}

\newcommand{\sobtg}[2]{\ensuremath{\mathcal{W}^{#1,#2}_g}}
\renewcommand{\leb}[1]{\ensuremath{\LL^{#1}}}
\newcommand{\dx}{\d \vec x}
\renewcommand{\ds}{\d s}
\renewcommand{\eenorm}[2]{\| {{#1}} \|_{\sob2{#2}_h(\W)}}
\renewcommand{\enorm}[2]{\| {{#1}} \|_{\sob1{#2}_h(\W)}}
\newcommand{\znorm}[2]{\| {{#1}} \|_{\leb{#2}_h(\W)}}
\newcommand{\rbar}{\overline R}
\makeatletter
\newcommand*{\rom}[1]{\text{\expandafter\@slowromancap\romannumeral #1@}}
\makeatother

\numberwithin{equation}{section}
\setboolean{showtodo}{true}
\setboolean{shownotes}{true}
\setboolean{showchanges}{false}
\setboolean{usemathrsfs}{false}
\usepackage{tkz-graph}
\usetikzlibrary{arrows}

\author{
  Nikos Katzourakis
}
\address{
  Nikos Katzourakis
  \thanks{
    Department of Mathematics and Statistics, Whiteknights, University of Reading, Reading RG6 6AX, UK
    {\tt{N.Katzourakis@reading.ac.uk}}.
}}

\author{
  Tristan Pryer
}
\address{
  Tristan Pryer
  \thanks{
    Department of Mathematics and Statistics, Whiteknights, University of Reading, Reading RG6 6AX, UK
    {\tt{T.Pryer@reading.ac.uk}}.
}}

\thanks{N.K. was partially supported through the EPSRC grant
  EP/N017412/1. T.P. was partially supported through the EPSRC grant
  EP/P000835/1.}
\title[FE approximation of $\infty$-BiHarmonic Functions]
      {On the numerical approximation of $p$-Biharmonic and $\infty$-Biharmonic functions}
\date{\today}

\begin{document}
\maketitle

\renewcommand{\thefootnote}{\fnsymbol{footnote}}
\renewcommand{\thefootnote}{\arabic{footnote}}

\begin{abstract}
  The $\infty$-Bilaplacian is a \emph{third order} fully nonlinear PDE given by
  \begin{equation*}
    \Delta^2_\infty  u\,
    :=
    (\Delta u)^3 | \D (\Delta u) |^2 = 0.
  \end{equation*}
  In this work we build a numerical method aimed at quantifying the
  nature of solutions to this problem which we call
  $\infty$-Biharmonic functions. For fixed $p$ we design a mixed
  finite element scheme for the pre-limiting equation, the
  $p$-Bilaplacian
  \begin{equation*}
    \Delta^2_p  u\,
    :=
    \Delta ( | \Delta u | ^{p-2} \Delta u ) = 0.
  \end{equation*}
  We prove convergence of the numerical solution to the weak solution
  of $\Delta^2_p u = 0$ and show that we are able to pass to the limit
  $p\to\infty$. We perform various tests aimed at understanding the
  nature of solutions of $\Delta^2_\infty u$ and we prove convergence
  of our discretisation to an appropriate weak solution concept of
  this problem, that of $\mathcal D$-solutions.
\end{abstract}

\section{Introduction and the $\infty$-Bilaplacian}
\label{sec:introduction}

Let $\W \subset \reals^d$ be an open and bounded set. For a given
function $u : \W \to \reals$ we denote the gradient of $u$ as
$\D u : \W \to \reals^d$ and its Hessian
$\Hess u:\W \to \reals^{d\times d}$ and Laplacian
$\Delta u : \W \to \reals$. The $p$--Bilaplacian
\begin{equation}
  \label{eq:p-lap}
  \Delta^2_p  u\,
  :=
  \Delta\qp{\norm{\Delta u}^{p-2} \Delta u} = 0
\end{equation}
is a fourth order elliptic partial
differential equation (PDE) which is a nonlinear generalisation of the
Bilaplacian. Such problems typically arise from areas of elasticity,
in particular, the nonlinear case can be used as a model for travelling
waves in suspension bridges
\cite{LazerMcKenna:1990,GyulovMorosanu:2010}. It is a fourth order
analogue to its second order sibling, the $p$--Laplacian, and as such it is
useful as a prototypical nonlinear fourth order problem.

The efficient numerical simulation of general fourth order problems
has attracted growing interest. A conforming approach to this class of
problems would require the use of $\cont{1}$ finite elements, the
Argyris element for example \cite[Section 6]{Ciarlet:1978}. From a
practical point of view the approach presents difficulties, in that
the $\cont{1}$ finite elements are difficult to design and complicated
to implement, especially when working in three spatial
dimensions. Other possibilities include discontinuous Galerkin
methods, which form a class of nonconforming finite element method. If
$p=2$ we have the special case that the ($2$--)Bilaplacian, $\Delta^2
u = 0$, is linear.  It has been well studied in the context of both
$\cont{1}$ finite elements \cite{Ciarlet:1978} and discontinuous
Galerkin methods; for example, the papers \cite{LasisSuli:2003,
  GeorgoulisHouston:2009} study the use of $h$--$k$ dG finite elements
(where $k$ here means the local polynomial degree as opposed to the
usual convention which is $p$) applied to the
($2$--)Bilaplacian. Alternative methods do exist, including those of
virtual element type
\cite{BC:2013,CangianiGeorgoulisPryerSutton:2017} and recovered
element type \cite{GeorgoulisPryer:2018}. In addition to this, the
classical work of \cite{babuska:1980} proposed mixed methods for the
linear problem whose analysis was based on the mesh-dependent norms in
\cite{babuska:1980a}. The numerical approximation of $p$-Bilaplacian
(quasi-linear, fourth order) type PDEs is relatively untouched. To the
authors' knowledge, the only known work is \cite{Pryer:2014} where a
discontinuous Galerkin method based on a variational principle was
derived and was shown to converge under minimal regularity. However,
no rates of convergence were proven.

In this work we propose a method based on $\cont{0}$-mixed finite
elements very much in the spirit of \cite{babuska:1980}. We rewrite
the minimisation problem in mixed formulation and prove that the
method converges under minimal regularity of the solution. In
addition, using an inf-sup condition inspired by
\cite{Makridakis:2018,GeorgoulisPryer:2017,GeorgoulisMakridakisPryer:2017}
and tools from \cite{Sandri:1993,Farhloul:1998,GiraultRaviart:2012},
we are able to show that under additional regularity assumptions the
approximation converges with specific rates that depend on $p$.

Making use of these convergence results and the uniqueness of
solutions in one dimension from \cite{KatzourakisPryer:2018b}, extended
to multi-spatial dimensions in \cite{KatzourakisMoser:2018} we are
able to justify that approximations of the $p$-Bilaplacian for large
$p$ are ``good'' approximations to $\infty$-Biharmonic
functions. These functions are solutions of the $\infty$-Bilaplacian
which is the PDE
\begin{equation}
  \label{eq:infbilap}
  \Delta^2_\infty  u\,
  :=
  \qp{\Delta u}^3 \norm{\D \qp{\Delta u}}^2 = 0,
\end{equation}
derived in \cite{KatzourakisPryer:2018b} as the formal limit of the
$p$-Bilaplacian (\ref{eq:p-lap}) as $p \rightarrow \infty$. The
$\infty$-Bilaplacian is the prototypical example of a PDE from second
order Calculus of Variations in $\leb{\infty}$, arising as the
analogue of the \emph{Euler--Lagrange} equation associated with
critical points of the supremal functional
\begin{equation}
  \cJ[u;\infty] := \Norm{\Delta u}_{\leb{\infty}(\W)}.
\end{equation}

Variational problems in $\leb{\infty}$ are notoriously
challenging. The first order case is reasonably well understood and
was initiated in the sequence of works by Aronsson starting with
\cite{Aronsson:1965}. In this case, the respective Euler-Lagrange
equation associated with critical points of the functional
\begin{equation}
  \cJ[u] = \Norm{\D u}_{\leb{\infty}(\W)},
\end{equation}
is quasi-linear, second order and given by
\begin{equation}
  \label{eq:inflap}
  \Delta_\infty u = \frob{\qp{\D u \otimes \D u}}{\Hess u} = 0.
\end{equation}
This equation is called the $\infty$-Laplacian and can be derived
through a $p$-approximation of the underlying $\sob{1}{p}$ energy functional,
see \cite{Pryer:2015,KatzourakisPryer:2016}.

It can be shown that solutions to \eqref{eq:infbilap} can not, in
general, be $C^3$ even when $d=1$; in particular, the Dirichlet
problem is not solvable in the class of classical solutions. For a
more extensive discussion we refer to
\cite{KatzourakisPryer:2018b}. Hence, the development of a solution
concept which can be interpreted in an appropriate weak sense is in
order. In the case of the $\infty$-Laplacian, the appropriate notion
is that of the Crandall-Ishii-Lions notion of viscosity solutions
\cite{CrandallIshiiLions:1992}. For an introduction to this theory we
refer to the monograph \cite{Katzourakis:2015}. We note that in the
framework of viscosity solutions we can obtain uniqueness of solution
for the Dirichlet problem \cite{Jensen:1993}. In the case of second
order Calculus of Variations in $\leb{\infty}$ the viscosity solution
concept for the resulting equations is no longer applicable since we
do not have access to a maximum principle for third order PDEs like
\eqref{eq:infbilap}, from which the solution concept stems.

One possibility for a generalised solution concept to
(\ref{eq:infbilap}) is that of $\mathcal D$-solutions
\cite{Katzourakis:2016a, Katzourakis:2016b, KatzourakisPryer:2018b}.  Roughly, this is a
probabilistic approach where derivatives that do not exist
classically are represented as limits of difference quotients into
Young measures over a compactification of the space of derivatives. This solution concept has already borne substantial fruit in the first order vectorial case of Calculus of Variations in $\leb{\infty}$, as well as for more general PDE systems. In the present second order setting it proves to be an appropriate notion as well, since absolute minimisers $u \in \mathcal{W}^{2,\infty}_g(\W)$ satisfying
\begin{equation}
\label{AM}
  \Norm{\Delta u}_{\leb{\infty}(\W')} \leq \Norm{\Delta v}_{\leb{\infty}(\W')}
  \Foreach \W' \Subset \W \text{ and } v\in \mathcal{W}^{2,\infty}_g(\W')  , 
\end{equation}
are indeed unique $\mathcal D$-solutions of (\ref{eq:infbilap}). Note
that the appropriate space to take minimisers is \emph{not}
$\sobg{2}{\infty}(\Omega)$ but rather the larger space
\begin{equation} \label{1.7}
  \sobtg{2}{\infty}(\W) := \ensemble{u \in \intersection{p\in (1, \infty)} \sobg{2}{p}(\W)}{\Delta u \in \leb{\infty}(\W)}.
\end{equation}
In \cite{KatzourakisPryer:2018b} it has been shown that in one spatial
dimension the problem does indeed have a unique absolutely minimising
$\mathcal D$-solution and in \cite{KatzourakisMoser:2018} for higher
spatial dimension.

The design of numerical schemes that are compatible with these
solution concepts that are inherently incompatible with duality
techniques is extremely difficult. Even for the well developed area of
viscosity solutions most numerical schemes that exist which are
compatible with the solution concept are based on the arguments of
\cite{BarlesSouganidis:1991} which advocates approximations based on
differences satisfying a discrete monotonicity property. The only
other methodology in the design of numerical schemes for the
$\infty$-Laplacian is to make use of the variational principle from
which the equation is derived. Galerkin approximations of the
$p$-Laplacian can then be shown to converge to the viscosity solution
of the $\infty$-Laplacian \cite{Pryer:2015}. This method has also been
used to characterise the nature of solutions to the variational
$\infty$-Laplace system \cite{KatzourakisPryer:2016}. This is also the
approach we use here. We build a scheme convergent to the weak
solution of the $p$-Bilaplacian and then justify its use as an
approximation of $\infty$-Biharmonic functions. This allows us
significant insight as to the nature of non-classical solutions of the
$\infty$-Bilaplacian and to make various conjectures about their
structure and behaviour.

The rest of the paper is set out as follows: In \S \ref{sec:plap} we
formalise notation and begin exploring some of the properties of the
$p$-Bilaplacian. In particular, we reformulate the PDE as a saddle
point type problem. We show inf-sup conditions for the underlying operators
guarantee that the saddle point type problem is well posed, motivating the
discretisation of this directly. In \S
\ref{sec:fem} we perform the discretisation for fixed $p$ and show that
discrete versions of the inf-sup conditions hold. A priori results for
both primal and auxiliary variables are a consequence of this.
Numerical experiments are given in \S \ref{sec:numerics} illustrating
the behaviour of numerical approximations to this problem. In addition,
we examine the solutions for large $p$ and make various conjectures as
to the structure of solutions in multiple spatial dimensions.

\section{Approximation via the $p$-Bilaplacian}
\label{sec:plap}

In this section we describe how $\infty$-Biharmonic functions can be
approximated using $p$-Biharmonic functions. We give a brief
introduction to the $p$--Bilaplacian problem, beginning by introducing
the Sobolev spaces
\begin{gather}
  \leb{p}(\W)
  =
  \ensemble{\phi \text{ measurable}}
           {\int_\W \norm{\phi}^p \dx < \infty}   \text{ for } p\in[1,\infty) 
             \text{ and }
             \\
             \leb{\infty}(\W)
             =
             \ensemble{\phi \text{ measurable}}
                      {\esssup_\W \norm{\phi} < \infty},
                      \\
                      \sob{l}{p}(\W) 
                      = 
                      \ensemble{\phi\in\leb{p}(\W)}
                               {\D^{\vec\alpha}\phi\in\leb{p}(\W), \text{ for } \norm{\geovec\alpha}\leq l}
                               \text{ and }
    \sobh{l}(\W)
    := 
    \sob{l}{2}(\W),
\end{gather}
which are equipped with the following norms and semi-norms:
\begin{gather}
  \Norm{v}_{\leb{p}(\W)}^p
  :=
       {\int_\W \norm{v}^p} \dx \text{ for } p \in [1,\infty)
         \text{ and }
  \Norm{v}_{\leb{\infty}(\W)}
  :=
  \esssup_\W |v|
       \\
       \Norm{v}_{\sob{l}{p}(\W)}^p
       := 
       \sum_{\norm{\vec \alpha}\leq l}\Norm{\D^{\vec \alpha} v}_{\leb{p}(\W)}^p
       \\
       \norm{v}_{\sob{l}{p}(\W)}^p
       :=
       \sum_{\norm{\vec \alpha} = l}\Norm{\D^{\vec \alpha} v}_{\leb{p}(\W)}^p
\end{gather}
where $\vec\alpha = \{ \alpha_1,\dots,\alpha_d\}$ is a
multi-index, $\norm{\vec\alpha} = \sum_{i=1}^d\alpha_i$ and
derivatives $\D^{\vec\alpha}$ are understood in the weak sense. We pay
particular attention to the case $l = 2$ and define
\begin{gather}
  \sobg{2}{p}(\W) 
  :=
  g + \sobz{2}{p}(\W)
  =
   \ensemble{\phi\in\sob{2}{p}(\W)}{\phi\vert_{\partial \W} = g \text{ and } \D\phi\vert_{\partial \W} = \D g},
\end{gather}
for a prescribed function $g\in\sob{2}{\infty}(\W)$, where the boundary condition is understood in the trace sense if $\partial\W \in \cont{0,1}(\W)$. We note that if $p>d$, then the boundary condition is satisfied in the pointwise sense since $\sobz{2}{p}(\W) \subseteq \mathrm{C}^{1}(\overline{\W})$. 

For the $p$--Bilaplacian, the action functional is given as
\begin{equation}
  \cJ[u;p]
  =
  \int_\W \norm{\Delta u}^p \dx. \footnote{Typically $\cJ[u;p] = \tfrac{1}{p} \int \norm{\Delta u}^p$. Note here the rescaling has no effect on the resultant Euler--Lagrange equations as the Lagrangian is independent of $u$.}
\end{equation}
We then look to find a
minimiser over the space $\sobg{2}{p}(\W)$, that is, to find
$u\in\sobg{2}{p}(\W)$ such that
\begin{equation}
  \cJ[u;p] = \min_{v\in\sobg{2}{p}(\W)} \cJ[v;p].
\end{equation}

If we assume temporarily that we have access to a smooth minimiser,
\ie $u\in\cont{4}(\W)$, then, given that the Lagrangian is of second order,
we have that the Euler--Lagrange equations are (in general) fourth
order and read
\begin{equation}
  \label{eq:p-biharm}
  \Delta\qp{\norm{\Delta u}^{p-2}\Delta u} = 0.
\end{equation}
Note that, for $p=2$, the PDE reduces to the Bilaplacian
$\Delta^2 u = 0$. In general, the Dirichlet problem for the $p$-Bilaplacian is,
given $g\in\sob{2}{\infty}(\W)$, to find $u$ such that
\begin{equation}
  \label{eq:plap}
  \left\{\ \ \
  \begin{split}
    \Delta_p u:= \Delta\qp{\norm{\Delta u}^{p-2} \Delta u} &= 0, \ \ \ \text{ in } \W,
    \\
    u &= g, \ \ \ \text{ on } \partial\W,
    \\
    \D u &= \D g, \ \text{ on } \partial\W.
  \end{split}
  \right.
\end{equation}

\begin{Defn}[weak solution]
  The problem (\ref{eq:plap}) has a weak formulation. Consider the semilinear form
  \begin{gather}
    \bi{u}{v} := \int_\W \qp{\norm{\Delta u}^{p-2}\Delta u} \Delta v \dx .
  \end{gather}
 Then, $u \in
  \sobg{2}{p}(\W)$ is a \emph{weak solution} of (\ref{eq:plap}) if it satisfies
  \begin{equation}
    \bi{u}{v} = 0  \Foreach v \in \sobz{2}{p}(\W).
  \end{equation}
\end{Defn}

\begin{Pro}[coercivity of $\cJ$]
  \label{prop:coercive}
  Suppose that $u\in\sobz{2}{p}(\W)$ and $f\in\leb{q}(\W)$, where
  $\tfrac{1}{p}+\tfrac{1}{q}=1$. We have that the action functional
  $\cJ[\ \cdot \ ;p]$ is coercive over $\sobz{2}{p}(\W)$, that is,
  \begin{equation}
    \cJ[u; p] 
    \geq 
    C
    \norm{u}^p_{2,p}
    -\gamma,
  \end{equation}
  for some $C >0 \AND \gamma \geq 0$. Equivalently, we have that there
  exists a constant $C>0$ such that
  \begin{equation}
    \label{eq:coercive-bilinear}
    \bi{v}{v} \geq
    C
    \norm{v}_{2,p}^p \Foreach v\in\sobz{2}{p}(\W).
  \end{equation}
\end{Pro}

\begin{Cor}[weak lower semicontinuity]
  \label{cor:semicont}
  The action functional $\cJ$ is weakly lower semi-continuous
  over $\sobg{2}{p}(\W)$. That is, given a sequence of functions $\{
  u_j \}_{j\in\naturals}$ which has a weak limit $u\in\sobg{2}{p}(\W)$, we have
  \begin{equation}
    \cJ[u; p] \leq \liminf_{j\to\infty} \cJ[u_j; p].
  \end{equation}
\end{Cor}
\begin{Proof}
  The proof of this fact is a straightforward extension of \cite[Section 8.2 Thm
    1]{Evans:1998} to second order Lagrangians, noting that
  $\cJ$ is coercive (from Proposition \ref{prop:coercive}) and
  convex. We omit the full details for brevity.
\end{Proof}

\begin{Cor}[existence and uniqueness]
  \label{cor:unique}
  There exists a unique minimiser to the $p$--Dirichlet energy
  functional. Equivalently, there exists a unique (weak) solution
  $u\in\sobg{2}{p}(\W)$ to the (weak form of the) Euler--Lagrange
  equations:
  \begin{equation}
    \bi{u}{v}
    =
    \int_\W \norm{\Delta u}^{p-2} \Delta u \Delta v \dx
    =
    0
    \Foreach v\in\sobz{2}{p}(\W).
  \end{equation}
\end{Cor}
\begin{Proof}
  Again, the result can be deduced by extending the arguments in
  \cite[Section 8.2]{Evans:1998} or \cite[Thm 5.3.1]{Ciarlet:1978},
  again, noting the results of Propositions \ref{prop:coercive} and
  convexity. The full argument is omitted for brevity.
\end{Proof}

\begin{The}[the limit as $p\to\infty$]
  \label{the:ptoinf}
  Let $(u_p)_1^\infty$ denote a sequence of weak solutions $u_p
  \in\sobg{2}{p}(\W)$ to the $p$-Bilaplacian. Then, there exists a
  subsequence converging uniformly together with their derivatives to
  a (candidate $\infty$-Biharmonic) function
  $u_\infty\in\sobtg{2}{\infty}(\W)$. Namely,
  \begin{equation}
    u_{p_j} \to u_\infty \text{ in } C^{1}(\overline{\Omega}),
  \end{equation}
  along a subsequence as $p\to \infty$.
\end{The}
\begin{Proof}
  Let $u_p \in \sobg{2}{p}(\W)$ denote the weak solution of
  (\ref{eq:plap}). In view of Corollary \ref{cor:unique}, we know that
  $u_p$ minimises the energy functional
  \begin{equation}
    \cJ[u_p] = \int_\W \norm{\Delta u_p}^p \d \vec x.
  \end{equation}
  In particular,
  \begin{equation}
    \cJ[u_p] \leq \cJ[g],
  \end{equation}
  where $g \in \sob{2}{\infty}(\W)$ is the associated boundary data to
  (\ref{eq:plap}). Using this fact, we have
  \begin{equation}
    \Norm{\Delta u_p}_{\leb{p}(\W)}^p
    =
    \cJ[u_p]
    \leq
    \cJ[g]
    =
    \Norm{\Delta g}_{\leb{p}(\W)}^p,
  \end{equation}
  and we may infer that
  \begin{equation}
    \label{eq:limitpf1}
    \Norm{\Delta u_p}_{\leb{p}(\W)}
    \leq
    \Norm{\Delta g}_{\leb{p}(\W)}.
  \end{equation}
  Now fix a $k > d$ and take $p \geq k$. Then, by using H\"older's inequality with $r = \tfrac{p}{k}$ and $q = \tfrac{r}{r-1}$ such that $\tfrac{1}{r} + \tfrac{1}{q} = 1$, we obtain
  \begin{equation}
    \label{eq:limitpf3}
    \Norm{\Delta u_p}_{\leb{k}(\W)}^k
    =
    \int_\W \norm{\Delta u_p}^k \d \vec x
    \leq
    \qp{\int_\W 1^q \d \vec x}^{1/q} \qp{\int_\W \norm{\Delta u_p}^p \d \vec x}^{1/r}.
  \end{equation}
 Hence
  \begin{equation}
    \Norm{\Delta u_p}_{\leb{k}(\W)}^k
    \leq
    \norm{\W}^{\tfrac{r-1}{r}} \Norm{\Delta u_p}^{k}_{\leb{p}(\W)}
    =
    \norm{\W}^{1-\tfrac{k}{p}} \Norm{\Delta u_p}^{k}_{\leb{p}(\W)}
  \end{equation}
  and we see
  \begin{equation}
    \label{eq:limitpf2}
    \Norm{\Delta u_p}_{\leb{k}(\W)}
    \leq
    \norm{\W}^{\tfrac 1 k-\tfrac{1}{p}} \Norm{\Delta u_p}_{\leb{p}(\W)}.    
  \end{equation}
  By using the triangle inequality, a double application of the Poincar\'e
  inequality (since both $u = g$ and $\D u = \D g$ on $\partial \W$) from Proposition \ref{pro:poincare} and
  the Calderon-Zygmund $L^k$ estimates from Proposition \ref{prop:equiv-of-norms}, we have 
  \begin{equation}
    \begin{split}
    \Norm{u_p}_{\leb{k}(\W)}
    &\leq
    \Norm{u_p - g}_{\leb{k}(\W)}
    +
    \Norm{g}_{\leb{k}(\W)}
    \\
    &\leq
    C'(k,\W) \Norm{\Hess u_p - \Hess g}_{\leb{k}(\W)}
    +
    \Norm{g}_{\leb{k}(\W)}
    \\
    &\leq
    C(k,\W) \Norm{\Delta u_p - \Delta g}_{\leb{k}(\W)}
    +
    \Norm{g}_{\leb{k}(\W)}.
    \end{split}
  \end{equation}
 By utilising the triangle inequality again, we have
  \begin{equation}
    \begin{split}
    \Norm{u_p}_{\leb{k}(\W)}
    &\leq
    C \qp{\Norm{\Delta u_p}_{\leb{k}(\W)}
    +
    \Norm{g}_{\sob{2}{k}(\W)}}
    \\
    &\leq
    C \qp{\norm{\W}^{\tfrac{1}{k}-\tfrac{1}{p}} \Norm{\Delta u_p}_{\leb{p}(\W)}
    +
    \Norm{g}_{\sob{2}{k}(\W)}},
    \end{split}
  \end{equation}
by virtue of (\ref{eq:limitpf2}). Similarly, one may show that
  \begin{equation}
    \Norm{\D u_p}_{\leb{k}(\W)}
    \leq
    C \qp{
      \norm{\W}^{\tfrac{1}{k}-\tfrac{1}{p}} \Norm{\Delta u_p}_{\leb{p}(\W)}
      +
      \Norm{g}_{\sob{2}{k}(\W)} }.
  \end{equation}
  Thus, in view of (\ref{eq:limitpf1}) we infer that
  \begin{equation}
    \label{eq:limitpf4}
    \begin{split}
      \Norm{u_p}_{\sob{2}{k}(\W)}
      &\leq
      C \Norm{g}_{\sob{2}{k}(\W)}.
    \end{split}
  \end{equation}
  This means that for any $k > d$ we have the uniform bound
  \begin{equation}
    \sup_{p > k} \Norm{u_p}_{\sob{2}{k}(\W)} \leq C = C(k,\W).
  \end{equation}
By invoking standard weak compactness arguments, we may extract a sub-sequence $\{
  u_{p_j} \}_{j=1}^\infty \subset \{ u_{p} \}_{p=1}^\infty$ and a
  function $u_\infty\in\sob{2}{k}(\W)$ such that, for any $k > n$,
  \begin{equation}
    u_{p_j} \rightharpoonup u_\infty \text{ weakly in } \sob{2}{k}(\W)
  \end{equation}
 as $j\to \infty$ and
  \begin{equation}
    \begin{split}
      \Norm{u_\infty}_{\sob{2}{k}(\W)} &\leq \liminf_{j\to\infty} \Norm{u_{p_j}}_{\sob{2}{k}(\W)}
      \\
      &\leq
      \liminf_{j\to\infty} 
      C\Norm{g}_{\sob{2}{k}(\W)}.
    \end{split}
  \end{equation}
  Since this is true for any fixed $k$, it is clear that
  $u_\infty \in \intersection{k\in (1,\infty)} \sob{2}{k}(\W)$.
  Further, by the weak lower semi-continuity of the $\leb{k}$ norm, from
  (\ref{eq:limitpf2}) we may infer
  $\Delta u_\infty \in \leb{\infty}(\W)$ and hence
  $u_\infty \in\sobtg{2}{\infty}(\W)$, therefore concluding the proof.
\end{Proof}

\begin{Rem}[elementary properties]
  \label{pro:pq}
  We will throughout this exposition use the notation $p$ to denote the exponent appearing in the Lagrangian and $q$ its conjugate exponent which satisfies 
  \begin{equation}
    \frac{1}{p} + \frac{1}{q} = 1.
  \end{equation}
  For a given $v\in\leb{p}(\W)$ it then holds that
  \begin{equation}
    \Norm{\norm{v}^{p-1}}_{\leb{q}(\W)} = \Norm{v}^{p-1}_{\leb{p}(\W)}.
  \end{equation}
\end{Rem}

\begin{Pro}[Poincar\'e inequality]
  \label{pro:poincare}
 Let $\W \subset \reals^d$ be a bounded domain. For any $p\in [1,\infty]$, there exists a constant $C = C(\W,p)>0$ depending only on $\W$ and $p$ such that
    \begin{equation}
      \Norm{u}_{\leb{p}(\W)} \leq C(\W,p) \Norm{\D u}_{\leb{p}(\W)},
    \end{equation}
    for all $u\in\sob{1}{p}_0(\W)$.
\end{Pro}

\begin{Pro}[Calderon-Zygmund estimate {\cite[Cor 9.10]{Gilbarg:1983}}]
\label{prop:equiv-of-norms}
Let $\W \subset \reals^d$ be a bounded domain. Then, for any $p\in (1,\infty)$, there is a constant $C=C(d,p)>0$ depending only on $d$ and $p$ such that
\begin{equation}
  \Norm{\Hess u}_{\leb{p}(\W)} \leq C(d,p) \Norm{\Delta u}_{\leb{p}(\W)},
\end{equation}
for all $u\in\sobz{2}{p}(\Omega)$.
\end{Pro}

An immediate consequence of Propositions \ref{pro:poincare} and \ref{prop:equiv-of-norms} above is that the norm $\Norm{\cdot}_{2,p}$ is equivalent to either of the seminorms $\Norm{\Hess ( \cdot )}_{\leb{p}(\W)}$ and $\Norm{\Delta ( \cdot )}_{\leb{p}(\W)}$ over the space $\sobz{2}{p}(\W)$.

\subsection{Mixed formulation of the $p$-Bilaplacian}

The mixed formulation we propose to analyse is based on the
observation that if $\phi(t) = |t|^{p-2}t$, the inverse is well defined
as $\phi^{-1}(t) = \mathrm{sgn}(t)|t|^{1/\qp{p-1}} = |t|^{q-2}t$. Using this we make the following choice
of auxiliary variable
\begin{equation}
  w = \norm{\Delta u}^{p-2} \Delta u
\end{equation}
from which we can infer that 
\begin{equation}
  \norm{w}^{q-2} w = \Delta u.
\end{equation}
This allows us to write the problem as the mixed system:
\begin{equation}
  \left\{
    \begin{split}
      -\Delta u &= \norm{w}^{q-2} w,
      \\
      -\Delta w &= 0.
    \end{split}
  \right.
\end{equation}
The mixed formulation can be written in {a strong form} as: Find a pair
{$\qp{u, w} \in \sobg{2}{p}(\W) \times \leb{q}(\W)$} such that
\begin{equation}
  \label{eq:saddle}
 \ \ \ \left\{ \ \
  \begin{split}
    a(w,\psi) + b(u,\psi) &= 0,
    \\
    b(\phi, w) &= 0,
   \ \ \  \Foreach
    {\qp{\psi, \phi} \in \leb{q}(\W) \times \sobz{2}{p}(\W)},
  \end{split}
  \right.
\end{equation}
where the semilinear form $a(w,\psi)$ and bilinear form $b(u,\psi)$ are
given by
\begin{equation}
  \label{eq:bilinearforms}
  \left\{
   \ \
  \begin{split}
    a(w,\psi) &:= \int_\W \norm{w}^{q-2} w \psi \dx
    \\
    b(u,\psi) &:= {\int_\W -\Delta u \psi \dx}
  \end{split}
  \right.
\end{equation}

Notice that the problem (\ref{eq:p-lap}) has been reformulated in a
mixed form. Although we already know that the problem has a unique
solution as a consequence of Corollary \ref{cor:unique}, we will show
that the equivalent mixed formulation also admits a unique solution
since the methodology will be useful henceforth. We begin with the
following result.

\begin{Pro}[{Inf-sup stability of $b(\cdot, \cdot)$ over $\sobz 2p(\W)$}]
  \label{lem:infsupb}
  {For any $u_0\in\sobz{2}{p}(\W)$, the bilinear form $b(\cdot, \cdot)$
  satisfies the following inf-sup property:
  \begin{equation}
    \Norm{\Delta u_0}_{\leb{p}(\W)} \leq
    C \sup_{0 \neq v \in \leb{q}(\W)} \frac{b(u_0, v)}{\Norm{v}_{\leb{q}(\W)}}.
  \end{equation}}  
\end{Pro}
\begin{Proof}
  {Fix $u_0\in\sobz 2p(\W)$. Then, we certainly have that
    $\norm{\Delta u_0}^{p-2}\Delta u_0 \in \leb{q}(\W)$. Therefore, by
    choosing $v = \norm{\Delta u_0}^{p-2}\Delta u_0$ we have
  \begin{equation}
    b(u_0,v) = \Norm{\Delta u}^p_{\leb{p}(\W)}
  \end{equation}
  and that
  \begin{equation}
    \Norm{v}_{\leb{q}(\W)}
    =
    \Norm{\Delta u_0^{p-1}}_{\leb{q}(\W)}
    =
    \Norm{\Delta u_0}^{p-1}_{\leb{p}(\W)},
  \end{equation}
  in view of the property given in Remark \ref{pro:pq}. Hence we have
  \begin{equation}
    b(u_0,v)
    =
    \Norm{\Delta u_0}^p_{\leb{p}(\W)}
    =
    \Norm{\Delta u_0}_{\leb{p}(\W)}
    \Norm{v}_{\leb{q}(\W)},
  \end{equation}
  which implies the desired result.}
\end{Proof}

\begin{The}[The mixed formulation is well posed]
  \label{the:saddle-well-posed}
  For every $g\in \sob{2}{\infty}(\W)$, there exists a unique pair
  $(u,w)$ solving (\ref{eq:saddle}) that satisfies
  \begin{equation}
    \Norm{\Delta u}_{\leb{p}(\W)}
    +
    \Norm{w}_{\leb{q}(\W)}^{q-1}
    \leq
    C{
      \Norm{\Delta g}_{\leb{p}(\W)}
    }.
  \end{equation}
\end{The}
\begin{Proof}
  The results of Proposition \ref{lem:infsupb} show that, for $u_0 := u
  - g \in \sobz{2}{p}(\W)$, we have
  \begin{equation}
    \begin{split}
      \Norm{\Delta u_0}_{\leb{p}(\W)}
      &\leq
      \sup_{0 \neq v \in \leb{q}(\W)} \frac{b(u_0, v)}{\Norm{v}_{\leb{q}(\W)}}
      \\
      &\leq
      \sup_{0 \neq v \in \leb{q}(\W)} \frac{b(u, v)}{\Norm{v}_{\leb{q}(\W)}}
      +
      \sup_{0 \neq v \in \leb{q}(\W)} \frac{b(g, v)}{\Norm{v}_{\leb{q}(\W)}}
      \\
      &\leq
      \sup_{0 \neq v \in \leb{q}(\W)} \frac{-a(w, v)}{\Norm{v}_{\leb{q}(\W)}}
      +
      \sup_{0 \neq v \in \leb{q}(\W)} \frac{b(g, v)}{\Norm{v}_{\leb{q}(\W)}}
      .
    \end{split}
  \end{equation}
  in view of (\ref{eq:saddle}). Now, by using Remark \ref{pro:pq} we
  estimate
  \begin{equation}
    \label{eq:ex-un1}
    \begin{split}
      \Norm{\Delta u_0}_{\leb{p}(\W)}
      &\leq
      \Norm{w^{q-1}}_{\leb{p}(\W)}
      +
      \Norm{\Delta g}_{\leb{p}(\W)}
      \\
      &\leq
      {\Norm{w}^{q-1}_{\leb{q}(\W)}
      +
      \Norm{\Delta g}_{\leb{p}(\W)}
      }.
    \end{split}
  \end{equation}  
  Now take $\psi = w$ in (\ref{eq:saddle}). Then,
  \begin{equation}
    \begin{split}
      a(w,w) + b(u,w) = 0.
    \end{split}
  \end{equation}
  Set $\phi = u_0$ in (\ref{eq:saddle}). Then,
  \begin{equation}
    b(u_0,w) = 0
  \end{equation}
  and in particular
  \begin{equation}
    \begin{split}
      a(w,w) + b(u,w) - b(u_0,w) = 0.
    \end{split}
  \end{equation}
  This in turn implies 
  \begin{equation}
    a(w,w) + b(g,w) = 0,
  \end{equation}
  or explicitly
  \begin{equation}
    \int_\W \norm{w}^q - \Delta g w \d \vec x = 0.
  \end{equation}
  Hence
  \begin{equation}
    \begin{split}
      \Norm{w}^q_{\leb{q}(\W)} &= \int_\W \Delta g w \d \vec x
      \\
      &\leq
      \Norm{\Delta g}_{\leb{p}(\W)} \Norm{w}_{\leb{q}(\W)},
    \end{split}
  \end{equation}
  and
  \begin{equation}
    \Norm{w}^{q-1}_{\leb{q}(\W)}
    \leq
    \Norm{\Delta g}_{\leb{p}(\W)},
  \end{equation}
  which yields the desired result upon noting
  \begin{equation}
    \Norm{\Delta u}_{\leb{p}(\W)}
    \leq
    \Norm{\Delta u_0}_{\leb{p}(\W)}
    +
    \Norm{\Delta g}_{\leb{p}(\W)}
  \end{equation}
  and combining with (\ref{eq:ex-un1}).
\end{Proof}

\begin{Rem}[Convergence to ``weak'' solutions to the $\infty$-Bilaplacian]
  Theorem \ref{the:ptoinf} guarantees convergence to a candidate
  $\infty$-Harmonic function. The correct notion of weak solution to
  the limiting problem
  \begin{equation} \label{2.59}
\ \  \left\{\ \ 
    \begin{split}
      \qp{\Delta u}^3 \norm{\D \qp{\Delta u}}^2
      &= 0, \ \ \text{ in }\W,
      \\
      u &= g, \ \ \text{ on } \partial\W,
      \\
      \D u &= \D g, \text{ on } \partial\W,
    \end{split}
    \right.
  \end{equation}
  is that of $\mathcal D$-solutions
  \cite{KatzourakisPryer:2018b,KatzourakisMoser:2018}. The solution is
  probabilistic in nature and interpreted in a weak sense. It is the
  only candidate $\infty$-Biharmonic function which means Theorem
  \ref{the:ptoinf} guarantees convergence of the sequence of
  $p$-Biharmonic functions to the unique $\infty$-Biharmonic $\mathcal
  D$-solution.
\end{Rem}

\section{Discretisation of the $p$-Bilaplacian}
\label{sec:fem}

In this section we describe a mixed finite element discretisation of the $p$-Bilaplacian.
Let $\T{}$ be a conforming triangulation of $\W$,
namely, $\T{}$ is a finite family of sets such that
\begin{enumerate}
\item $K\in\T{}$ implies $K$ is an open simplex (segment for $d=1$,
  triangle for $d=2$, tetrahedron for $d=3$),
\item for any $K,J\in\T{}$ we have that $\closure K\meet\closure J$ is
  a full lower-dimensional simplex (i.e., it is either $\emptyset$, a vertex, an
  edge, a face, or the whole of $\closure K$ and $\closure J$),
\item $\union{K\in\T{}}\closure K=\closure\W$.
\end{enumerate}
The shape regularity constant of $\T{}$ is defined as the number
\begin{equation}
  \label{eqn:def:shape-regularity}
  \mu(\T{}) := \inf_{K\in\T{}} \frac{\rho_K}{h_K},
\end{equation}
where $\rho_K$ is the radius of the largest ball contained inside
$K$ and $h_K$ is the diameter of $K$. An indexed family of
triangulations $\setof{\T n}_n$ is called \emph{shape regular} if 
\begin{equation}
  \label{eqn:def:family-shape-regularity}
  \mu:=\inf_n\mu(\T n)>0.
\end{equation}
We let $\E{}$ be the skeleton (set of common interfaces) of the
triangulation $\T{}$ and say $e\in\E$ if $e$ is on the interior of
$\W$ and $e\in\partial\W$ if $e$ lies on the boundary $\partial\W$.

We let $\poly k(\T{})$ denote the space of piecewise polynomials of
degree $k\geq 2$ over the triangulation $\T{}$, that is,
\begin{equation}
  \poly k (\T{}) = \{ \phi \text{ such that } \phi|_K \in \poly k (K) \},
\end{equation}
 and introduce the \emph{finite element space}
\begin{gather}
  \label{eqn:def:finite-element-space}
  \fes := \poly k(\T{}) \cap \cont{0}(\W),
\end{gather}
to be the usual space of continuous piecewise polynomial
functions. {We define jump operators for arbitrary scalar
  functions $v$ and vectors $\vec v$ over an edge $e$ shared by
  elements $K_1$ and $K_2$ as $\jump{v} = {{{v}|_{K_1} \geovec n_{K_1}
      + {v}|_{K_2}} \geovec n_{K_2}}$, $\jump{\vec v} =
  {{\vec{v}|_{K_1}}} \cdot \geovec n_{K_1} + {{\vec{v}|_{K_2}}} \cdot
  \geovec n_{K_2}$ and when $e$ is on $\partial \W$ we understand
  $\jump{v} = {v}|_{K} \geovec n_{\partial\W}$ and $\jump{\vec v} =
          {\vec v}|_{K}\cdot \geovec n_{\partial\W}$.}

Further, we define $\funk h\W\reals$ to be the {piecewise
  constant} \emph{meshsize function} of $\T{}$ given by
\begin{equation}
  h(\vec{x}):=\max_{\closure K\ni \vec{x}}h_K.
\end{equation}
A mesh is called quasi-uniform when there exists a positive constant
$C$ such that $\max_{x\in\Omega} h \le C \min_{x\in\Omega} h$. In what
follows we shall assume that all triangulations are shape-regular and
quasi-uniform although the results may be extendable even in the
non-quasi-uniform case using techniques developed in
\cite{DieningKreuzer:2008}.

\begin{Defn}[Ritz projection operators]
  \label{def:ritz}
  The Ritz projection operator $R$ is defined through requiring
  \begin{equation}
    \int_\W \D \qp{R v} \cdot \D \phi \d \vec x
    =
    \int_\W \D v \cdot \D \phi \d \vec x
    \Foreach \phi\in\fes\cap \sobh1_0(\W),
  \end{equation}
  and $Rv$ coincides with an appropriate interpolant of $v$ on the
  boundary. This operator satisfies the following approximation
  properties for quasi-uniform meshes \cite{Li:2017}: for any $v \in
  \sob{k+1}{q}(\W)$, and $k\geq2$
  \begin{gather}
    \Norm{v - R v}_{\leb{q}(\W)}
    +
    \Norm{h\qp{\D v - \D \qp{ R v}}}_{\leb{q}(\W)}
    +
    \qp{\sum_{K\in\T{}} \Norm{h^2\qp{\Delta v - \Delta \qp{ R v}}}_{\leb{q}(K)}^q}^{1/q}
    \leq C h^{k+1} \norm{v}_{k+1,q}.
  \end{gather}
  The Neumann Ritz projection $\rbar$ is defined through requiring
  orthogonality over a larger space
  \begin{equation}
    \int_\W \D \qp{\rbar w} \cdot \D \psi \dx = \int_\W \D w \cdot \D \psi \dx \Foreach \psi\in\fes
  \end{equation}
  and requiring
  \begin{equation}
    \int_\W \rbar w \dx = \int_\W w \dx.
  \end{equation}
  The results of \cite{Li:2017} also imply that $\rbar$ satisfies the
  same approximation properties as $R$.
\end{Defn}

\begin{Defn}[Mesh-dependent norms]
  {We introduce the {mesh-dependent} $\leb{p}$- and $\sob 2 p$-norms to be
    \begin{equation}
      \begin{split}
        \znorm{w_h}{p}^p
        &:=
        \Norm{w_h}_{\leb{p}(\W)}^p
        +
        \Norm{h^{1/p}w_h}_{\leb{p}(\E)}^p 
        \\
        \eenorm{w_h}{p}^p
        &:=
        \Norm{\Delta_h w_h}_{\leb{p}(\W)}^p
        +
        \Norm{h^{1/p-1}\jump{\D w_h}}_{\leb{p}(\E)}^p,
      \end{split}
    \end{equation}
    where $\Delta_h$ denotes an elementwise Laplace operator. }
\end{Defn}

\subsection{Galerkin discretisation}
  Consider the space
  \begin{equation}
    \fes_g := \{ \phi \in \fes : \phi \vert_{\partial \Omega} = R g \}.
  \end{equation}
  Then, we consider the Galerkin discretisation of (\ref{eq:plap}), to find $\qp{u_h, w_h} \in \fes_g \times \fes$ such that
\begin{equation}
  \label{eq:plapdis}
  \begin{split}
    a(w_h,\psi) + b_h(u_h,\psi) &= 0
    \\
    b_h(\phi, w_h) &= 0, \ \ \ \Foreach \qp{\psi, \phi} \in \fes \times \fes_0,
  \end{split}
\end{equation}
where the bilinear form $a(\cdot, \cdot)$ is given in
(\ref{eq:bilinearforms}), $b_h(\cdot, \cdot)$ is a consistent
discretisation of $b(\cdot, \cdot)$ given by {
  \begin{equation}
    b_h(u_h, \psi)
    =
    -\sum_{K\in\T{}} \int_K \Delta u_h \psi \dx
    +
    \int_\E \jump{\D u_h} \psi \ds.
  \end{equation}
  Notice that the method is equivalent to finding $\qp{u_h, w_h} \in
  \fes_g \times \fes$ such that
  \begin{equation}
    \begin{split}
      \int_\W \norm{w_h}^{q-2} w_h \psi + \D u_h \cdot \D \psi \dx
      &=
      \int_{\partial \W} \D g \cdot \vec n \psi \ds
      \\
      \int_\W \D w_h \cdot \D \phi \dx& = 0
      , \ \ \ \Foreach \qp{\psi, \phi} \in \fes \times \fes_0.
    \end{split}
  \end{equation}
  Hence the Ritz projection operator from Definition \ref{def:ritz} is
  the $b_h$- orthogonal projection onto $\fes_g$, that is, for
  $v\in\sobh1_g(\W)$
  \begin{equation}
    b_h(Rv - v, \phi) = 0 \Foreach \phi \in \fes_0.
  \end{equation}
}

\begin{Rem}
  \label{rem:boundedness}
  {
    The reason for defining the mesh-dependent norms as we do is to
    ensure the boundedness property
  \begin{equation}
    \norm{b_h(u_h, v_h)}
    \leq
    \eenorm{u_h}{p} \znorm{v_h}{q}.
  \end{equation}
  The scaling in the edge terms is chosen so that for arbitrary $v_h
  \in \fes$ each mesh-dependent norm is equivalent to the continuous
  counterpart, that is $\znorm{v_h}{p} \sim \Norm{v_h}_{\leb{p}(\W)}$
  for example.}
\end{Rem}

\begin{Lem}
  \label{lem:infsupbdis}
  { Assume the mesh is quasi-uniform, then the bilinear form
    $b_h(\cdot, \cdot)$ satisfies the following inf-sup property: for
    any $\Phi \in \fes_0$,
  \begin{equation}
    \eenorm{\Phi}{p}
    \leq
    C
    \sup_{0 \neq v_h \in \fes_0}
    \frac{b_h(\Phi, v_h)}
         {\znorm{v_h}{q}}.
  \end{equation}   }
\end{Lem}
\begin{Proof} 
  { The proof of this fact takes inspiration from
    \cite{Makridakis:2018} (see also \cite{GeorgoulisPryer:2017} and
    \cite{GeorgoulisMakridakisPryer:2017} for related ideas). We begin
    by showing that there exists a function $v$ that is discrete but
    not an element of $\fes_0$ such that
  \begin{equation}
    \label{eq:step1}
    b_h(\Phi, R v) \geq C \eenorm{\Phi}{p}^p
  \end{equation}
  and then showing the discrete stability estimate that $\znorm{R
    v}{q}\leq C \eenorm{\Phi}{p}^{p-1}$ .}

  {
    To begin we denote $b_K$ as the cubic a posteriori bubble
  function. This is a function that is $\poly 3$, positive over $K$,
  extended by zero outside of $K$ and satisfies that
  $\Norm{b_K}_{\leb{\infty}(K)} = 1$. Now take $v_1$ such that $v_1|_K
  = -b_K \norm{\Delta \Phi}^{p-2} \Delta \Phi$. Notice that
  $v_1\in\sobz 1 q(\W)$ and that $v_1|_e = 0$ for all $e\in\E \cup \partial \W$. Then
  through the equivalence of norms over finite dimensional linear
  spaces.
  \begin{equation}
    \label{eq:infsup0}
    \begin{split}
      \frac 1 C
      \sum_{K\in\T{}} \Norm{\Delta\Phi}^p_{\leb{p}(K)}
      \leq
      \sum_{K\in\T{}} \int_K b_K \norm{\Delta
        \Phi}^p \dx
      =
      b_h(\Phi, v_1)
      =
      b_h(\Phi, R v_1).
    \end{split}
  \end{equation}}

  {Now let $b_e$ be the edge bubble function that vanishes over all
  vertices of $\T{}$. Again this is a polynomial that is positive over
  $K$, extended by zero outside of the two elements sharing $e\in\E$
  and satisfies $\Norm{b_e}_{\leb{\infty}(e)} = 1$. Define $v_e : e
  \to \reals$ such that $v_e = h^{1-p} \norm{\jump{\D
      \Phi}}^{p-2}\jump{\D \Phi}$ on the face $e$ and extended by a
  constant on the direction normal to $e$. Set $v_2 := \sum_{e\in\E}
  b_e v_e$ then we have $v_2\in\sobz{1}{q}(\W)$ and
  \begin{equation}
    \begin{split}
      b_h(\Phi, v_2)
      &=
      \sum_{K\in\T{}} \int_K -\Delta \Phi v_2 \dx + \int_\E \jump{\D \Phi} v_2 \ds
      \\
      &=
      \sum_{K\in\T{}} \int_K -\Delta \Phi v_2 \dx + \int_\E b_e h^{1-p} \norm{\jump{\D \Phi}}^p \ds.
    \end{split}
  \end{equation}
  Now equivalence of norms shows there exists a constant $C>0$
  independent of $\Phi$ and $h$ such that
  \begin{equation}
    \label{eq:infsup1}
    \begin{split}
      \frac 1 C \Norm{h^{1/p-1} \jump{\D \Phi}}^p_{\leb{p}(\E)} 
      &\leq
      \int_\E b_e h^{1-p} \norm{\jump{\D \Phi}}^p \ds
      \\
      &=
      b_h(\Phi, v_2)
      +
      \sum_{K\in\T{}} \int_K \Delta \Phi v_2 \dx 
      \\
      &\leq
      b_h(\Phi, R v_2)
      +
      \qp{\sum_{K\in\T{}} \int_K \norm{\Delta \Phi}^p \dx}^{1/p} \Norm{v_2}_{\leb{q}(\W)}.
    \end{split}
  \end{equation}
  Young's inequality with $\epsilon$ shows that
  \begin{equation}
    \label{eq:infsup2}
    \begin{split}
      \qp{\sum_{K\in\T{}} \int_K \norm{\Delta \Phi}^p \dx}^{1/p} \Norm{v_2}_{\leb{q}(\W)}
      &\leq
      C(\epsilon) \qp{\sum_{K\in\T{}} \int_K \norm{\Delta \Phi}^p \dx}
      +
      \epsilon \Norm{v_2}_{\leb{q}(\W)}^q
      \\
      &\leq
      C(\epsilon) \qp{\sum_{K\in\T{}} \int_K \norm{\Delta \Phi}^p \dx}
      +
      C\epsilon \Norm{h^{1/p-1}\jump{\D \Phi}}_{\leb{p}(\E)}^p
    \end{split}
  \end{equation}
  in view of the definition of $v_2$. Now substituting
  (\ref{eq:infsup2}) into (\ref{eq:infsup1}) and choosing $\epsilon$
  appropriately small we see that
  \begin{equation}
    \label{eq:infsup3}
    \begin{split}
      \Norm{h^{1/p-1} \jump{\D \Phi}}^p_{\leb{p}(\E)} 
      &\leq
      C\qp{
      b_h(\Phi, R v_2)
      +
      \qp{\sum_{K\in\T{}} \int_K \norm{\Delta \Phi}^p \dx}}
      \\
      &\leq
      C\qp{
      b_h(\Phi, R v_2)
      +
      b_h(\Phi, R v_1)
      },
    \end{split}
  \end{equation}
  by (\ref{eq:infsup0}). Hence with $v=v_1+v_2$ we have shown
  (\ref{eq:step1}).}

  {We must now show the stability bound. To begin we show a stability
  result for the Ritz projection. With $z\in\sob{2}{p}(\W) \cap
  \sobh1_0(\W)$ solving the problem
  \begin{equation}
    -\Delta z = \norm{Rv}^{q-2} Rv,
  \end{equation}
  we see
  \begin{equation}
    \Norm{R v}_{\leb{q}(\W)}^q = b_h(z, Rv) = b_h(z, Rv - v) + b_h(z,
    v) = b_h(z - z_h, Rv - v) + b_h(z,v),
  \end{equation}
  with $z_h$ chosen as the Cl\'ement interpolant of $z$. Now using the
  definition of $z$ and approximation properties of $z_h$
  \begin{equation}
    \begin{split}
      \Norm{R v}_{\leb{q}(\W)}^q
      &\leq
      C\Norm{h \qp{\D (Rv) - \D v}}_{\leb{q}(\W)} \Norm{z}_{\sob{2}{p}(\W)}
      +
      \int_\W \norm{Rv}^{q-2} Rv v \dx
      \\
      &\leq
      C\Norm{Rv}^{q-1}_{\leb{q}(\W)}\qp{
        \Norm{h \qp{\D (Rv) - \D v}}_{\leb{q}(\W)}
      +
      \Norm{v}_{\leb{q}(\W)}}.
    \end{split}
  \end{equation}
  Using the $\sob{1}{q}$ stability of $R$ from \cite{Li:2017} we have
  \begin{equation}
    \Norm{R v}_{\leb{q}(\W)}
    \leq
    C\qp{\Norm{h\D v}_{\leb{q}(\W)}
    +
    \Norm{v}_{\leb{q}(\W)}}.
  \end{equation}
  Notice we have not used the super-approximation ideas from
  \cite{Makridakis:2018,GeorgoulisMakridakisPryer:2017} and are
  working on quasi-uniform meshes only. Now for $v=v_1+v_2$ defined
  above we are able to use inverse inequalities to see that
  \begin{equation}
    \Norm{R v}_{\leb{q}(\W)}
    \leq
    C(p) \Norm{v}_{\leb{q}(\W)},
  \end{equation}
  and through the definition of $v_1$ and $v_2$ we have
  \begin{equation}
    \Norm{v}_{\leb{q}(\W)}
    \leq
    C\eenorm{\Phi}{p}^{p-1}.
  \end{equation}
  Hence
  \begin{equation}
    \znorm{Rv}{q} \eenorm{\Phi}{p}
    \leq
    C \Norm{R v}_{\leb{q}(\W)} \eenorm{\Phi}{p}
    \leq
    C \Norm{v}_{\leb{q}(\W)} \eenorm{\Phi}{p}
    \leq
    C \eenorm{\Phi}{p}^p
    \leq
    C b_h(\Phi, Rv)
  \end{equation}
  and certainly
  \begin{equation}
    \eenorm{\Phi}{p}
    \leq
    C\frac {b_h(\Phi, Rv)}{\znorm{Rv}{q}}
    \leq
    \sup_{v_h\in\fes_0}
    C\frac {b_h(\Phi, v_h)}{\znorm{v_h}{q}},
  \end{equation}
  concluding the proof.}
\end{Proof}

\begin{The}[existence and uniqueness of solution to (\ref{eq:plapdis})]
  \label{the:exist-unique-dis}
  There exists a unique pair $\ \qp{u_h, w_h} \in \fes_g \times \fes$
  solving (\ref{eq:plapdis}). They satisfy the stability bound
  \begin{equation}
    \label{eq:dis-stab-bound}
    \eenorm{u_h}{p} + \Norm{w_h}^{q-1}_{\leb{q}(\W)}
    \leq
    C {
      \Norm{\Delta g}_{\leb{p}(\W)}
    }
    .
  \end{equation}
  Note that since $g\in\sob{2}{\infty}(\W)$, the right hand side of
  (\ref{eq:dis-stab-bound}) is finite.
\end{The}
\begin{Proof}
  The proof of this mirrors that of Theorem
  \ref{the:saddle-well-posed}. We begin by noting that for $\psi =
  w_h$ we have
  \begin{equation}
    \begin{split}
      a(w_h, w_h) + b_h(u_h, w_h) = 0.
    \end{split}
  \end{equation}
  Now for $\phi = u_{h,0} := u_h - R g$ we see that
  \begin{equation}
    b_h(u_h - R g, w_h) = 0, 
  \end{equation}
  hence
  \begin{equation}
    \begin{split}
      0
      &=
      a(w_h, w_h) + b_h(R g, w_h).
    \end{split}
  \end{equation}
  Now, by definition, we obtain
  \begin{equation}
    \label{pof:1}
    \begin{split}
      \Norm{w_h}^q_{\leb{q}(\W)}
      &\leq
      \eenorm{Rg}{p} \znorm{w_h}{q}
      \\
      &\leq
      C
      \Norm{\Delta  g}_{\leb{p}(\W)}
      \Norm{w_h}_{\leb{q}(\W)}
      ,
    \end{split}
  \end{equation}
  by Remark \ref{rem:boundedness} and Lemma \ref{lem:infsupbdis} and
  hence
  \begin{equation}
    \Norm{w_h}^{q-1}_{\leb{q}(\W)}
    \leq
    \Norm{\Delta  g}_{\leb{p}(\W)}.
  \end{equation}
  The result follows because
  \begin{equation}
    \label{pof:2}
    \begin{split}
      \eenorm{u_{h,0}}{p}
      &\leq
      C \sup_{0\neq v_h \in\fes_0} \frac{b_h(u_{h,0}, v_h)}{\znorm{v_h}{q}}
      \\
      &\leq
      C \qp{
        \sup_{0\neq v_h \in\fes_0} \frac{b_h(u_{h}, v_h)}{\znorm{v_h}{q}}
        +
        \sup_{0\neq v_h \in\fes_0} \frac{b_h(Rg, v_h)}{\znorm{v_h}{q}}       
      }
      \\
      &\leq
      C \qp{
        \sup_{0\neq v_h \in\fes_0} \frac{-a(w_h, v_h)}{\znorm{v_h}{q}}
        +
        \sup_{0\neq v_h \in\fes_0} \frac{b_h(g, v_h)}{\znorm{v_h}{q}}       
      }
      \\
      &\leq
      C \qp{\Norm{\norm{w_h}^{q-1}}_{\leb{p}(\W)} + \Norm{\Delta g}_{\leb{p}(\W)}}
      \\
      &\leq
      C \qp{\Norm{{w_h}}^{q-1}_{\leb{q}(\W)} + \Norm{\Delta g}_{\leb{p}(\W)}}
    \end{split}
  \end{equation}
  by the discrete inf-sup condition in Lemma \ref{lem:infsupbdis} and the same argument as in the proof of Theorem \ref{the:saddle-well-posed}. Since
  \begin{equation}
    \label{pof:3}
    \begin{split}
      \eenorm{u_{h}}{p} &\leq {\eenorm{u_{h,0}}{p} + \eenorm{Rg}{p}}
      \\
      &\leq C\qp{\eenorm{u_{h,0}}{p} + \Norm{\Delta g}_{\leb{p}(\W)}}
    \end{split}
  \end{equation}
  combining (\ref{pof:1}), (\ref{pof:2}) and (\ref{pof:3}) concludes the proof.
\end{Proof}

Next we state some technical properties that will be used in the theorem that follows.

\begin{Lem}[{Properties of $a(\cdot, \cdot)$, cf.\ \cite[Prop 3.1]{Sandri:1993}}]
  \label{lem:sandri}
  With $w\in\leb{q}(\W)$ and $w_h,v_h\in\fes$, for any $p \geq 2$,
  there exist constants
  \begin{enumerate}
  \item  $C_1>0$ such that
    \begin{equation}
      C_1 \frac{\Norm{w - w_h}_{\leb{q}(\W)}^2}
      {\Norm{w}_{\leb{q}(\W)}^{2-q} + \Norm{w_h}_{\leb{q}(\W)}^{2-q}}
      \leq
      a(w, w - w_h) - a(w_h,w - w_h).
    \end{equation}
  \item $C_2>0$ such that
    \begin{equation}
      C_2 \int_\W \norm{ \norm{w}^{q-2} w - \norm{w_h}^{q-2} w_h } \norm{w - w_h} \dx
      \leq
      a(w, w - w_h) - a(w_h,w - w_h).
    \end{equation}
  \item $C_3>0$ such that
    \begin{equation}
      a(w, w - v_h) - a(w_h, w - v_h)
      \leq
      C_3
      \qp{
        \int_\W \norm{\norm{w}^{q-2} w - \norm{w_h}^{q-2} w_h}
        \norm{w - w_h}\dx
      }^{1/p}
      \Norm{w - v_h}_{\leb{q}(\W)}.
    \end{equation}
  \end{enumerate}
\end{Lem}

\begin{The}[Approximability of the numerical schemes]
  Let $\qp{u, w}\in\sob{k+1}{p}_g(\W) \times \sob{k+1}{q}(\W)$ be the
  unique solution of (\ref{eq:saddle}) and $\qp{u_h, w_h} \in
  \fes_g\times \fes$ be the finite element approximation satisfying
  (\ref{eq:plapdis}). Then, the following error estimate holds
  \begin{equation}
    \begin{split}
      \Norm{w - w_h}_{\leb{q}(\W)}
      +
      \eenorm{u - u_h}{p}^{p-1}
      &\leq
      C \qp{
        h^{\frac q 2 \qp{k+1}} \norm{w}^{q/2}_{\sob{k+1}{q}(\W)}
        +
        h^{{k+1}} \norm{w}_{\sob{k+1}{q}(\W)}
        +
        h^{k-1} \norm{u}_{\sob{k+1}{p}(\W)}
      }.
    \end{split}
  \end{equation}
\end{The}

\begin{Proof}
  We begin by noting the Galerkin orthogonality results
  \begin{equation}
    \label{eq:Galerkinorthog}
    \begin{split}
      b_h(\phi, w - w_h) &= 0 \Foreach \phi\in{\fes_0},
      \\
      a(w, \psi) - a(w_h, \psi) + b_h(u - u_h, \psi) &= 0 \Foreach \psi\in\fes,
    \end{split}
  \end{equation}
  in view of (\ref{eq:saddle}) and (\ref{eq:plapdis}). 

  Now using Lemma \ref{lem:sandri} we have 
  \begin{equation}
    \label{eq:approx-proof-0}
    \begin{split}
      \frac{C_1\Norm{w - w_h}_{\leb{q}(\W)}^2}
           {2\qp{\Norm{w}_{\leb{q}(\W)}^{2-q} + \Norm{w_h}_{\leb{q}(\W)}^{2-q}}}
           +
      \frac{C_2}{2}
      \int_\W \norm{ \norm{w}^{q-2} w - \norm{w_h}^{q-2} w_h } \norm{w - w_h} \d \vec x
      &\leq
      a(w, w - w_h) - a(w_h,w - w_h)
      \end{split}
  \end{equation}
  Now using the semilinearity of $a(\cdot,\cdot)$ we have, for
  $\chi\in\fes$ denoting some approximation of $w$ to be chosen, that
  \begin{equation}
    \begin{split}
      a(w, w - w_h) - a(w_h,w - w_h)
      &=
      a(w, w - \chi) - a(w_h, w - \chi) + a(w, \chi - w_h) - a(w_h, \chi - w_h)
      \\
      &=
      \underbrace{a(w, w - \chi) - a(w_h, w - \chi)}_{=:\rom{1}}
      +
      \underbrace{b_h(u - u_h, w_h - \chi)}_{=:\rom{2}},
    \end{split}
  \end{equation}
  in view of (\ref{eq:Galerkinorthog}). We proceed to bound these
  terms separately, starting with $\rom{1}$.

  Making use of Lemma \ref{lem:sandri}
  \begin{equation}
    \label{eq:approx-proof-1}
    a(w, w-\chi) - a(w_h, w-\chi)
    \leq
    C_3 \qp{\int_\W \norm{\norm{w}^{q-2} w - \norm{w_h}^{q-2} w_h} \norm{w-w_h} \d \vec x}^{1/p}
    \Norm{w - \chi}_{\leb{q}(\W)}.
  \end{equation}
  Young's inequality with $\epsilon$ states for $a,b,\epsilon > 0$
  \begin{equation}
    ab \leq \frac 1 p \qp{\epsilon a}^p + \frac 1 q \qp{\frac b \epsilon}^q,
  \end{equation}
  which, upon applying to (\ref{eq:approx-proof-1}), shows
  \begin{equation}
    \label{eq:approx-proof-2}
    a(w, w-\chi) - a(\chi, w-\chi)
    \leq
    \frac{\epsilon^p}{p} 
    {\int_\W \norm{\norm{w}^{q-2} w - \norm{w_h}^{q-2} w_h} \norm{w-w_h} \d \vec x}
    +
    \frac{C_3^q}{q \epsilon^q}\Norm{w - \chi}_{\leb{q}(\W)}^q.
  \end{equation}
  Now choosing $\epsilon = \qp{\tfrac{C_2 p}{2}}^{1/p}$ and
  we have
  \begin{equation}
    \label{eq:approx-proof-3}
    a(w, w-\chi) - a(\chi, w-\chi)
    \leq
    \frac{C_2}{2} 
    {\int_\W \norm{\norm{w}^{q-2} w - \norm{w_h}^{q-2} w_h} \norm{w-w_h} \d \vec x}
    +
    C(q) \Norm{w - \chi}_{\leb{q}(\W)}^q.
  \end{equation}
  Notice we have picked $\epsilon$ such that the first term on the
  right hand side of (\ref{eq:approx-proof-3}) will cancel with the
  second term on the left hand side of (\ref{eq:approx-proof-0}).
  
  To control $\rom{2}$ we pick $\chi$ such that
  \begin{equation}
    b_h(\phi, \chi) = 0 \Foreach \phi\in\fes_0.
  \end{equation}
  An example of such an operator is the Neumann Ritz projection
  operator, $\rbar w$, given in Definition \ref{def:ritz}. With this
  choice of $\chi$, noting the definition of $w_h$ from
  (\ref{eq:plapdis}), it is clear that
  \begin{equation}
    b_h(\phi, w_h - \chi) = 0 \Foreach \phi\in\fes_0,
  \end{equation}
  and hence
  \begin{equation}
    b_h(u - u_h, w_h - \chi) = b_h(u-u_h - R\qp{u - u_h}, w_h - \chi) = b_h(u - Ru, w_h - \chi).
  \end{equation}
  Now making use of the boundedness of $b_h(\cdot, \cdot)$ we have
  \begin{equation}
    \label{eq:approx-proof-4}
    \begin{split}
      b_h(u - u_h, w_h - \chi)
      &\leq
      \eenorm{u-R u}{p} \znorm{w_h - \chi}{q}
      \\
      &\leq
      C
      \eenorm{u-R u}{p} \Norm{w_h - \chi}_{\leb{q}(\W)}
      \\
      &\leq
      \frac{C}{4 \epsilon }\eenorm{u-R u}{p}^2
      +
      \epsilon \Norm{w_h - \chi}_{\leb{q}(\W)}^2
      \\
      &\leq
      \frac{C}{4 \epsilon }\eenorm{u-Ru}{p}^2
      +
      2\epsilon \qp{
        \Norm{w - w_h}_{\leb{q}(\W)}^2
        + \Norm{w - \chi}_{\leb{q}(\W)}^2
      }.
    \end{split}
  \end{equation}
  Substituting (\ref{eq:approx-proof-3}) and (\ref{eq:approx-proof-4})
  into (\ref{eq:approx-proof-0}) and choosing $\epsilon$ small enough
  we see
  \begin{equation}
    \label{eq:approx-proof-5}
    \begin{split}
      \Norm{w - w_h}_{\leb{q}(\W)}^2
      &\leq
      C\qp{
        \Norm{w - \chi}_{\leb{q}(\W)}^q
        +
        \eenorm{u-Ru}{p}^2
        +
        \Norm{w - \chi}_{\leb{q}(\W)}^2
      },
    \end{split}
  \end{equation}
  allowing us to use the approximability of $R$ and $\rbar$ concluding
  the proof of the auxiliary variable.
  
  To show a bound for the primal variable we make use of the inf-sup
  condition from Lemma \ref{lem:infsupbdis}, noting that in view of
  Galerkin orthogonality and the definition of $R$ we have
  \begin{equation}
    \begin{split}
      0 &=
      a(w,\phi) - a(w_h,\phi) + b_h(u-u_h, \phi)
      \\
      &=
      a(w,\phi) - a(w_h,\phi) + b_h(R u-u_h, \phi)       \Foreach \phi\in\fes_0.
    \end{split}
  \end{equation}
  It is then clear that 
  \begin{equation}
    \label{eq:aid3}
    \begin{split}
      \eenorm{{R u} - u_h}{p} 
      &\leq 
      \sup_{0 \neq \phi \in \fes_0} 
      \frac {b_h\qp{Ru - u_h, \phi}}
      {\znorm{\phi}{q}}
      \\
      &=
      \sup_{0 \neq \phi \in \fes_0} 
      \frac {a\qp{w_h, \phi} - a\qp{w, \phi}}
      {\znorm{\phi}{q}}
      \\
      &\leq 
      C_3 \sup_{0 \neq \phi \in \fes_0} 
      \frac
      {\qp{\int_\W \norm{\norm{w}^{p-2} w - \norm{w_h}^{p-2} w_h} \norm{w - w_h} \dx}^{1/p} \Norm{\phi}_{\leb{q}(\W)}}
      {\znorm{\phi}{q}}
      \\
      &\leq 
      C_3 C \qp{\int_\W \norm{\norm{w}^{p-2} w - \norm{w_h}^{p-2} w_h} \norm{w - w_h} \dx}^{1/p},
    \end{split}
  \end{equation}
  through the equivalence of the $\leb{q}$-norm and its discrete
  counterpart. Now by Lemma \ref{lem:sandri} and Young's inequality with $\epsilon$ we have
  \begin{equation}
    \begin{split}
      C_2 {\int_\W \norm{\norm{w}^{p-2} w - \norm{w_h}^{p-2} w_h} \norm{w - w_h} \dx}
      &\leq
      a(w, w-w_h) - a(w_h, w-w_h)
      \\
      &\leq
      C_3 \qp{\int_\W \norm{\norm{w}^{p-2} w - \norm{w_h}^{p-2} w_h} \norm{w - w_h} \dx}^{1/p}
      \Norm{w-w_h}_{\leb{q}(\W)}
      \\
      &\leq
      \frac{\epsilon^p}{p} {\int_\W \norm{\norm{w}^{p-2} w - \norm{w_h}^{p-2} w_h} \norm{w - w_h} \dx}
      +
      \frac{C_3^q}{q \epsilon^q} \Norm{w-w_h}_{\leb{q}(\W)}^q.
    \end{split}
  \end{equation}
  The particular choice $\epsilon = \qp{\frac{pC_2}{2}}^{1/p}$ then shows that
  \begin{equation}
    \label{eq:aid2}
    {\int_\W \norm{\norm{w}^{p-2} w - \norm{w_h}^{p-2} w_h} \norm{w - w_h} \dx}
    \leq
    C \Norm{w-w_h}_{\leb{q}(\W)}^q.
  \end{equation}
  
  Substituting (\ref{eq:aid2}) into (\ref{eq:aid3}) results in
  \begin{equation}
    \eenorm{{R u} - u_h}{p} 
    \leq
    C\Norm{w - w_h}_{\leb{q}(\W)}^{q/p}.
  \end{equation}
  The result follows from the fact
  \begin{equation}
    \eenorm{u - u_h}{p} 
    \leq
    \eenorm{R u - u_h}{p} 
    +
    \eenorm{R u - u}{p} 
  \end{equation}
  and using the approximation properties of the Ritz projection, concluding the proof.
\end{Proof}

\begin{Rem}[Optimality of the bounds]
  Notice that the rates trail off as $p$ gets large. A similar
  phenomena was noticed when constructing methods for the
  $p$-Laplacian \cite[Thm 5.3.5]{Ciarlet:1978} where for a conforming
  piecewise linear approximation, $u_h$, the error behaved like
  \begin{equation}
    \Norm{u - u_h}_{\sob{1}{p}(\W)} \leq C h^{1/\qp{p-1}}.
  \end{equation}
  An analysis based on quasi-norms \cite{BarrettLiu:1994a} was then
  introduced to rectify this. It may be possible to use these
  techniques to show optimal error bounds for the $p$-Bilaplacian
  based on the quasi-norm
  \begin{equation}
    \Norm{u}_{v,p}^p
    :=
    \int_\W \norm{\Delta u}^2 \qp{\norm{\Delta u} + \norm{\Delta v}}^{p-2} \dx.
  \end{equation}
  We shall not push this point further in this work however. Instead,
  in order to try to characterise the limiting problem, we shall focus
  on convergence under minimal regularity.
\end{Rem}

We begin by defining the semilinear form
\begin{equation}
  c\qp{\qp{u,w},\qp{\phi,\psi}} := a(w,\psi) + b_h(u,\psi) + b_h(\phi,w),
\end{equation}
then the discrete mixed form of the Bilaplacian can be written, equivalently to
(\ref{eq:plapdis}), as seeking
$\qp{u_h,w_h} \in \fes_g \times \fes$ such that
\begin{equation}
  \label{eq:c}
  c\qp{\qp{u_h,w_h},\qp{\phi,\psi}} = 0 \Foreach \qp{\phi,\psi} \in \fes_0 \times \fes.
\end{equation}

\begin{The}[Convergence under minimal regularity]
  \label{the:convergenceundermin}
  Let $\qp{u_h,w_h}$ be a sequence of finite element solutions of (\ref{eq:plapdis}) indexed by the mesh parameter $h$ and let also $u\in\sob{2}{p}_g(\W)$ be the solution of the $p$-Bilaplacian. Then we have
  \begin{itemize}
  \item $u_h \to u \text{ strongly in }\leb{p}$ as $h\to 0$,
  \item $w_h \rightharpoonup w \text{ weakly in } \leb{q}$ as $h\to 0$.
  \end{itemize}
\end{The}
\begin{Proof}
  The stability result given in Theorem \ref{the:exist-unique-dis} allows us to infer that
  the sequence $\qp{u_h,w_h}$ is bounded uniformly in $h$. This means, up to a subsequence, that there exists a
  $\qp{u^*, w^*} \in \sobg{2}{p}(\W) \times \leb{q}(\W)$ such
  that $u_h \to u^*$ strongly in $\leb{p}(\W)$ and
  $w_h \rightharpoonup w^*$ weakly in $\leb{q}(\W)$.

  Now suppose $v_1\in\cont{\infty}(\W)$. Take
  $\qp{\phi,\psi} = \qp{0, \rbar v_1}$ in (\ref{eq:c}). Then,
  \begin{equation}
    0
    =
    c\qp{\qp{u_h,w_h},\qp{0,\rbar v_1}}
    =
    a(w_h, \rbar v_1) + b_h (u_h, \rbar v_1).
  \end{equation}
  Since $u_h \to u^*$ and by the properties of the projection $\rbar$
  given in Definition \ref{def:ritz} we have that
  \begin{equation}
    b_h(u_h, R v_1) \to b(u^*, v_1).
  \end{equation}
  Also, since $w_h \rightharpoonup w^*$ and $\rbar v_1 \to v_1$ strongly we have
  \begin{equation}
    a(w_h, R v_1) \to a(w^*, v_1).
  \end{equation}
  Hence
  \begin{equation}
    a(w_h, R v_1) + b_h(u_h, R v_1) \to  a(w^*, v_1) + b(u^*, v_1).
  \end{equation}
  Now suppose $v_2\in\cont{\infty}_0(\W)$ and take $\qp{\phi,\psi} = \qp{R v_2, 0}$ in (\ref{eq:c}), then
  \begin{equation}
    b_h(w_h, R v_2) = 0.
  \end{equation}
  By the same arguments we have
  \begin{equation}
    b_h (R v_2,w_h)
    \to 
    b(v_2,w^*).
  \end{equation}
  Using density of $\cont{\infty}_0(\W)\times\cont{\infty}(\W)$
  functions in $\sobz{2}{p}(\W) \times \leb{q}(\W)$ shows that
  $\qp{u^*, w^*}$ must solve the Bilaplacian and since the solution
  was unique, the whole sequence $\qp{u_h, w_h} \to \qp{u, w}$.
\end{Proof}

\begin{Cor}
  Let $u_{h,p} \in \fes_g$ be the Galerkin solution of
  (\ref{eq:plapdis}) and let $u_\infty$ denote a candidate
  $\infty$-Biharmonic function. Then, along a subsequence we have
  \begin{equation}
    u_{h,p_j} \to u_\infty \in \cont{0}(\overline{\Omega}) \text{ as } p\to\infty \text{ and } h\to 0. 
  \end{equation}
\end{Cor}

\begin{Rem}
  Since there exists a unique subsequential $p$-Biharmonic limit
  $u_\infty$ to the $\infty$-Bilaplacian on $\Omega$ the whole
  sequence must converge to this function, that is
  \begin{equation}
    u_{h,p} \to u_\infty \in \cont{0}(\overline{\Omega}) \text{ as } p\to\infty \text{ and } h\to 0. 
  \end{equation}
\end{Rem}

\section{Numerical experiments}
\label{sec:numerics}

In this section we summarise numerical experiments validating the analysis done in previous sections.

\subsection{Test 1: Benchmarking a $2$-dimensional problem.}

We begin by benchmarking the scheme against a known solution of the
$p$-Biharmonic problem. To do this we introduce a source term into the
problem
\begin{equation}
  \label{eq:plapsource}
  \left\{
  \ \ 
  \begin{split}
    \Delta\qp{\norm{\Delta u}^{p-2} \Delta u} &= f,  \ \ \text{ in } \W,
    \\
    u &= g, \ \ \text{ on } \partial\W,
    \\
    \D u &= \D g, \text{ on } \partial\W.
  \end{split}
  \right.
\end{equation}
This allows us to pick a function $g$ and construct the appropriate
source term such that $g$ solves (\ref{eq:plapsource}). For these
tests we choose
\begin{equation}
  u(x,y) = \frac1 {\pi^2}\sin{\pi x}\sin{\pi y}.
\end{equation}
We take $\W = [-1,1]^2$ and discretise the domain with a sequence of
concurrently refined criss-cross type meshes.

\begin{figure}[!ht]
  \caption[]
          {\label{fig:aronson}
            {\bf Test 1:} Benchmarking results for the mixed finite element
            approximation to (\ref{eq:plapsource}). We test the cases
            $p=2,\dots, 7$ for polynomials of degree $k=2$. The
            results show that the convergence rates as predicted in
            the analysis are achieved for the primal variable. Note
            that in the case $p > 2$ convergence rates are both higher
            than predicted for both primal \emph{and} auxiliary
            variable. Notice also that as $p$ increases the auxiliary
            variable converges at a faster rate.
  }
  \begin{center}
    \subfigure[{\label{fig:a2}
        The $2$-Bilaplacian.
    }]{
      \includegraphics[scale=\figscale,width=0.45\figwidth]{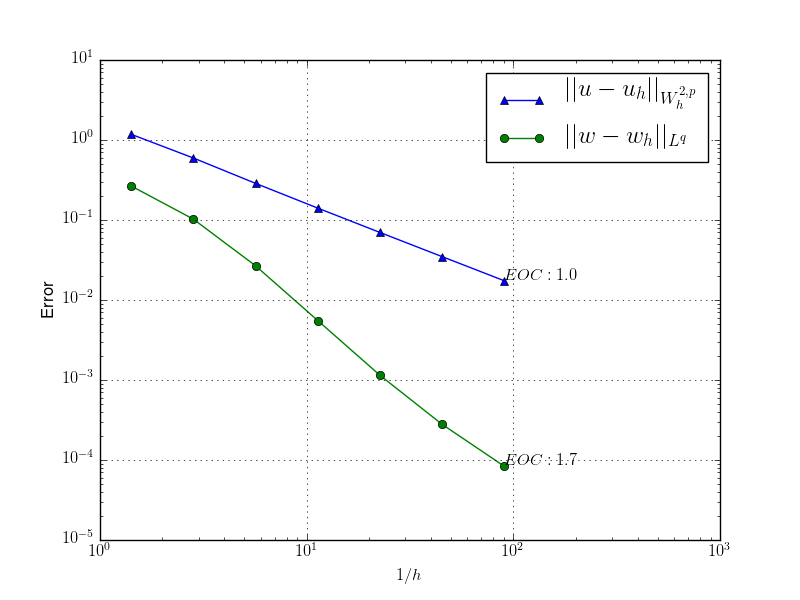}
    }
    \hfill
    \subfigure[{\label{fig:a1}
        The $3$-Bilaplacian.
    }]{
      \includegraphics[scale=\figscale,width=0.45\figwidth]{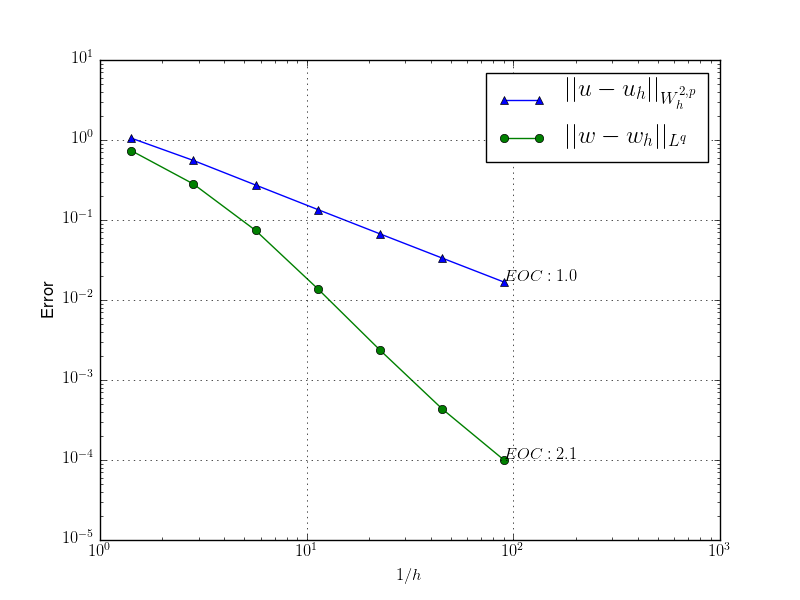}
    }
    \hfill
    \subfigure[{\label{fig:a2}
        The $4$-Bilaplacian.
    }]{
      \includegraphics[scale=\figscale,width=0.45\figwidth]{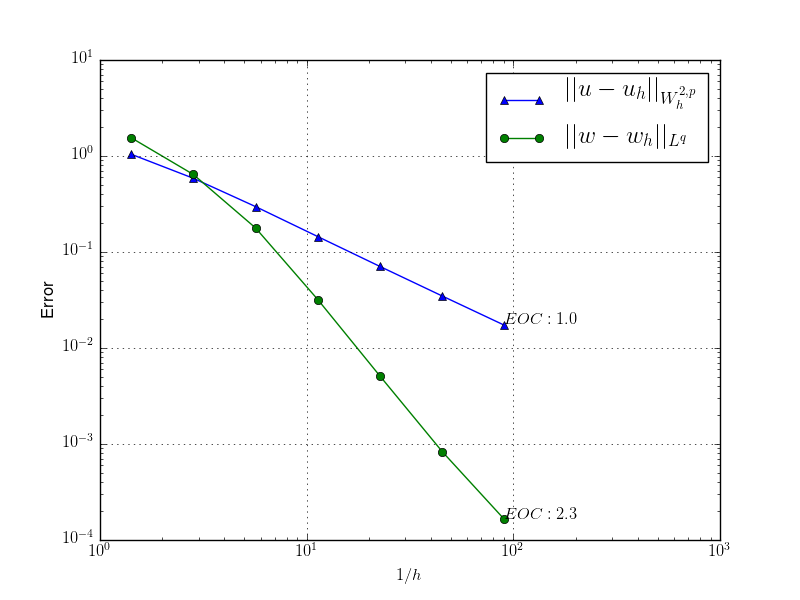}
    }
    \hfill
    \subfigure[{\label{fig:a2}
        The $5$-Bilaplacian.
    }]{
      \includegraphics[scale=\figscale,width=0.45\figwidth]{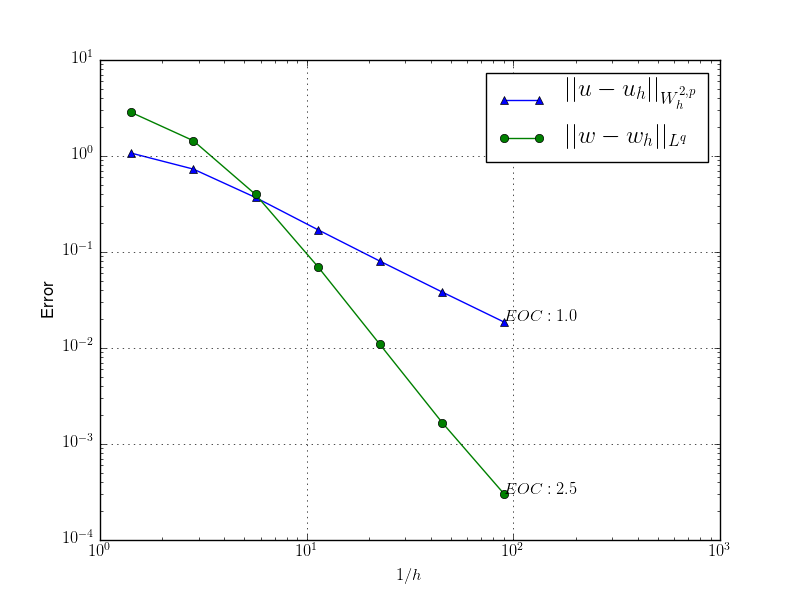}
    }
    \hfill
    \subfigure[{\label{fig:a2}
        The $6$-Bilaplacian.
    }]{
      \includegraphics[scale=\figscale,width=0.45\figwidth]{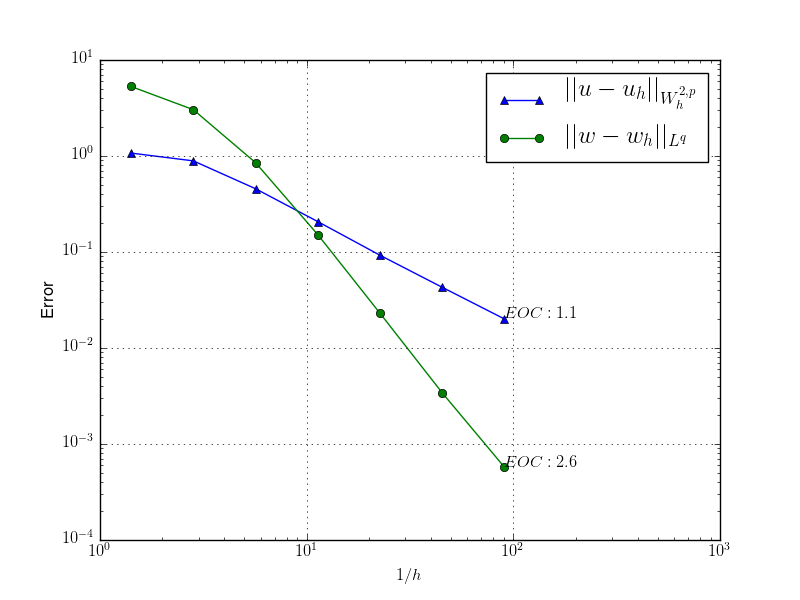}
    }
    \hfill
    \subfigure[{\label{fig:a2}
        The $7$-Bilaplacian.
    }]{
      \includegraphics[scale=\figscale,width=0.45\figwidth]{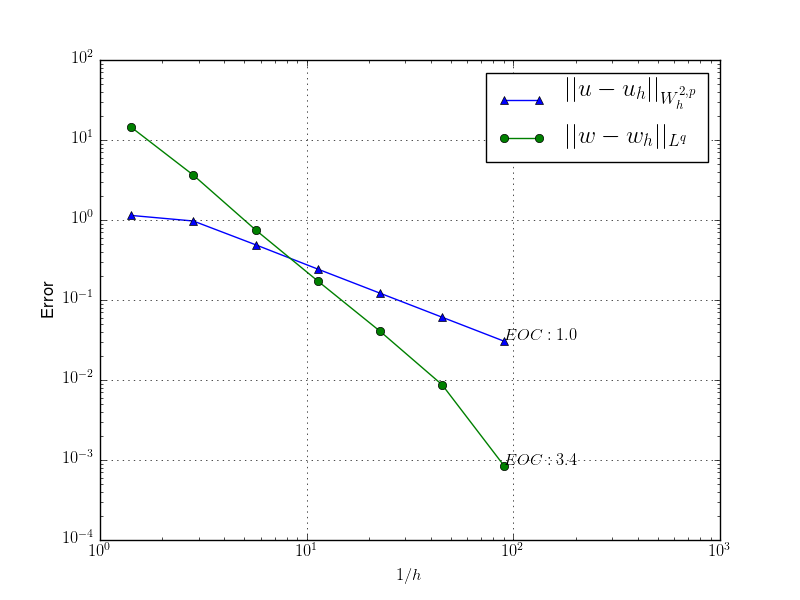}
    }
    \end{center}
  \end{figure}

The nonlinear system of equations generated are solved using a damped
Newton method initialised by solving the $2$-Bilapacian with
corresponding boundary data and forcing. The damping parameter is
chosen as $\tfrac 1 {p-2}$. The results are presented in Figure
\ref{fig:aronson}.

\subsection{Test 2: Characterising $\infty$-Harmonic functions in $1$-dimension}

In this experiment we illustrate some of the properties of
$\infty$-Biharmonic functions. The results illustrate that for
practical purposes, as one would expect, the approximation of
$p$-Biharmonic functions for large $p$ gives good resolution of
candidate $\infty$-Biharmonic functions.

We consider the Dirichlet problem for the $p$-Bilaplacian for $d=1$
with the boundary data given by the values of the cubic function
\begin{equation} \label{9.2}
  g(x) = \tfrac{1}{120} (4x-3)(2x-1)(4x-1)
\end{equation}
on $[0,1]$. We simulate the $p$-Bilaplacian for increasing values of
$p$ and present the results in Figure \ref{fig:floppydonkeydick}
indicating that in the limit the $\infty$-Biharmonic function should
be piecewise quadratic.

{
\begin{figure}[!ht]
  \caption[]
  {\label{fig:floppydonkeydick}
    {\bf Test 2:} A mixed finite element approximations to an $\infty$-Biharmonic function using $p$-Biharmonic functions for various $p$ for the problem given by \eqref{9.2}. Notice that as $p$ increases, $u''$ tends to a piecewise constant up to Gibbs oscillations. This is an indication the solution is indeed piecewise quadratic. Also there is only one breaking point in the solution, the location and size of this discontinuity was fully characterised in \cite{KatzourakisPryer:2018b}.
  }
  \begin{center}
    \subfigure[{\label{fig:a1}
        The approximation to $u$, the solution of the $4$-Bilaplacian.
    }]{
      \includegraphics[scale=\figscale,width=0.47\figwidth]{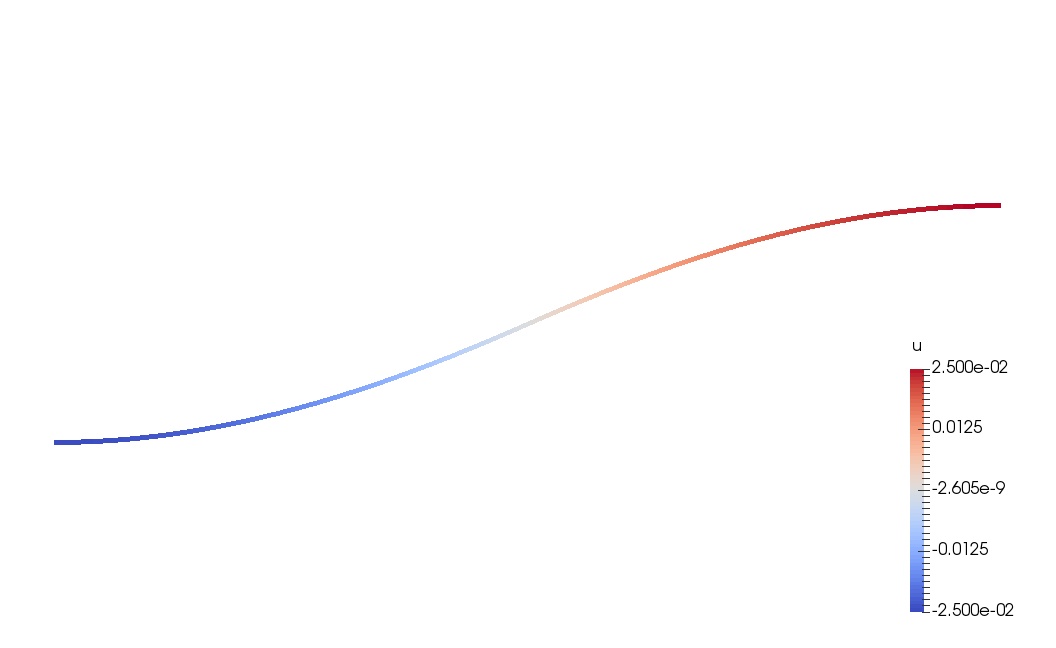}
    }
    \hfill
    \subfigure[{\label{fig:a2}
                The approximation to $u$, the solution of the $202$-Bilaplacian.
    }]{
      \includegraphics[scale=\figscale,width=0.4\figwidth]{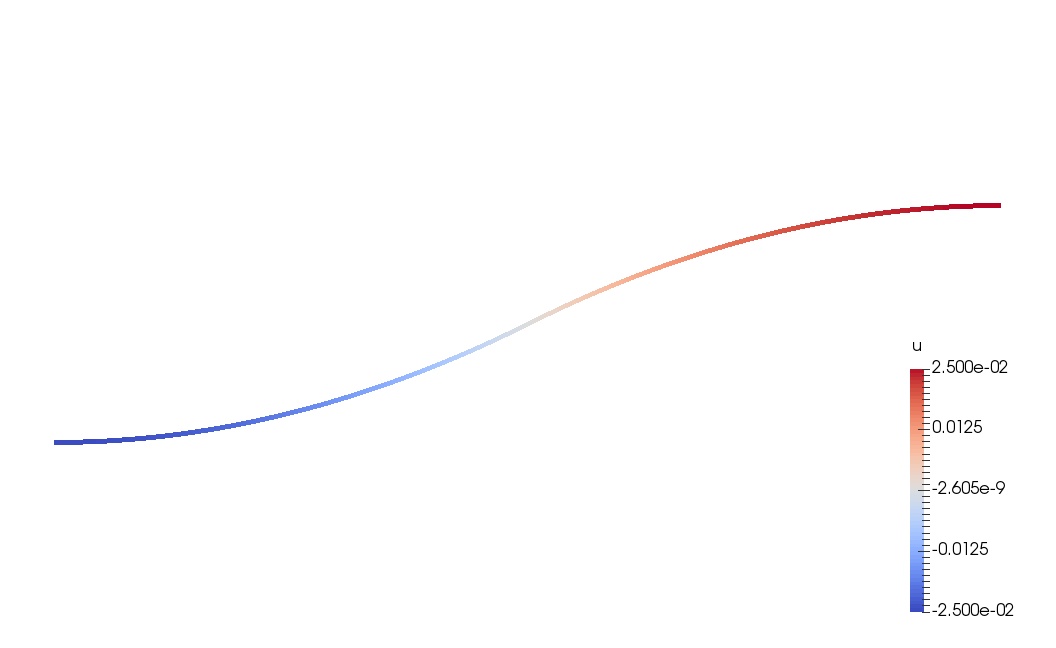}
    }    
   \subfigure[{\label{fig:a1}
        The approximation to $u''$, the Laplacian of the solution of the $4$-Bilaplacian.
    }]{
      \includegraphics[scale=\figscale,width=0.47\figwidth]{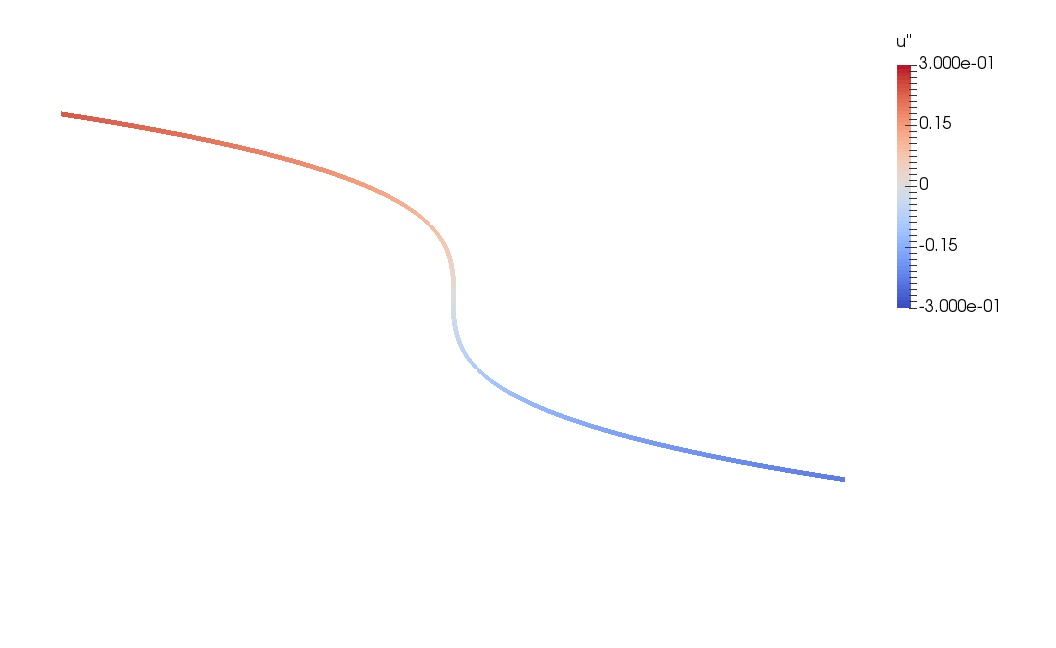}
    }
    \hfill
    \subfigure[{\label{fig:a2}
                The approximation to $u''$, the Laplacian of the solution of the $12$-Bilaplacian.
    }]{
      \includegraphics[scale=\figscale,width=0.47\figwidth]{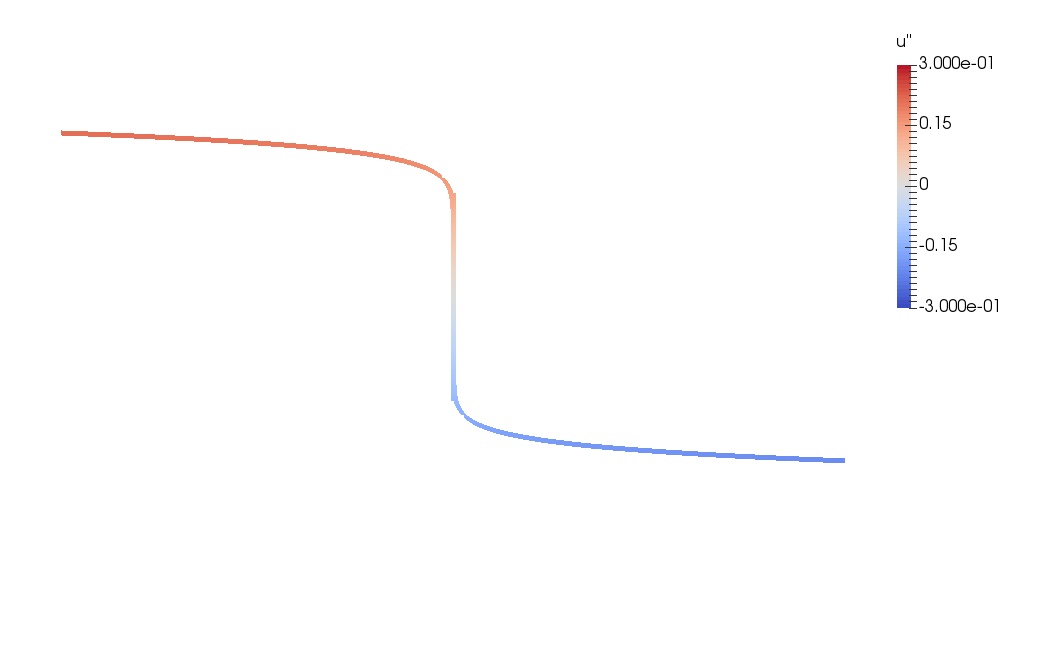}
    }
    \hfill
    \subfigure[{\label{fig:a2}
                The approximation to $u''$, the Laplacian of the solution of the $42$-Bilaplacian.        
    }]{
      \includegraphics[scale=\figscale,width=0.47\figwidth]{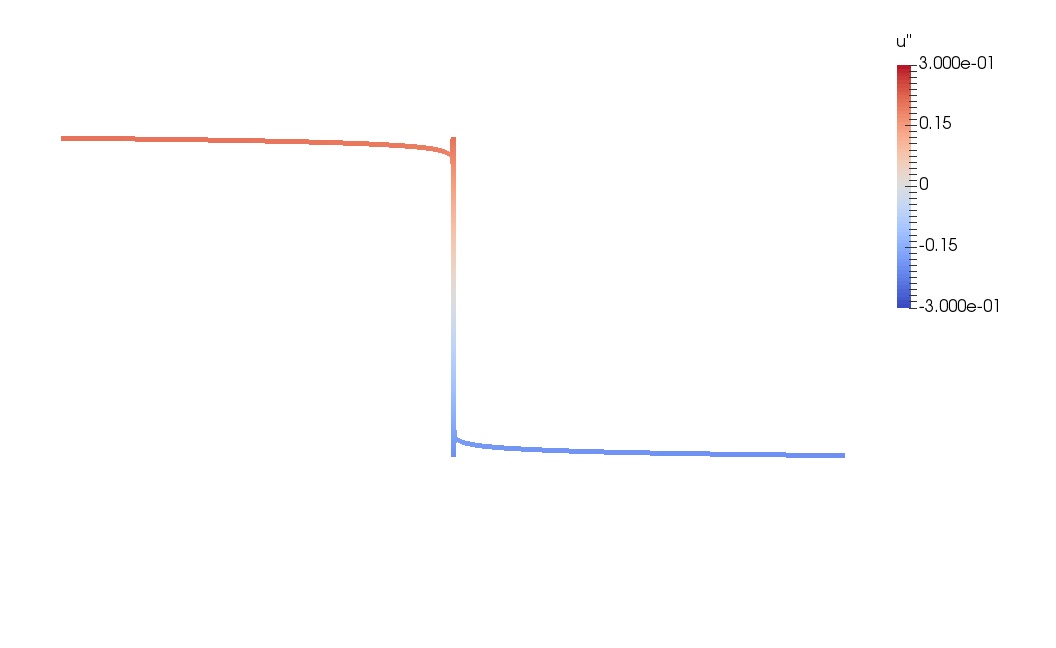}
    }
    \hfill
    \subfigure[{\label{fig:a2}
                        The approximation to $u''$, the Laplacian of the solution of the $202$-Bilaplacian.
    }]{
      \includegraphics[scale=\figscale,width=0.47\figwidth]{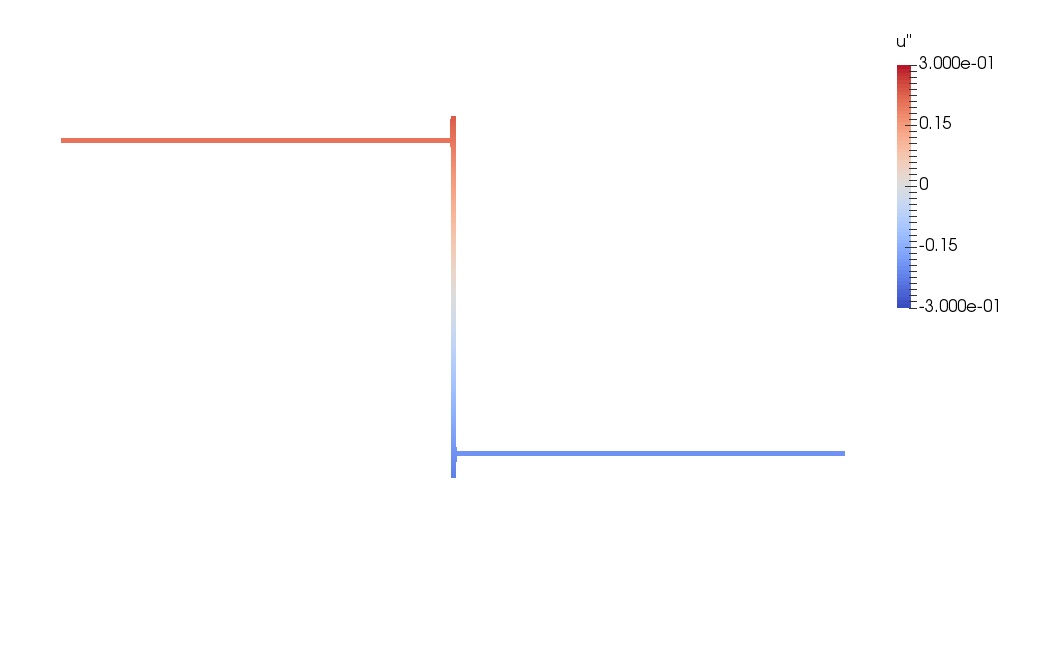}
    }
  \end{center}
  \end{figure}

}

\subsection{Test 3: Characterising $\infty$-Harmonic functions in $2$-dimensions}

\begin{figure}[!ht]
  \caption[]
  {\label{fig:floppydonkeydick4}
    {\bf Test 3b:} A mixed finite element approximations to an $\infty$-Biharmonic function using $p$-Biharmonic functions for various $p$ for the problem given by \eqref{eq:2dpbiharm} and \eqref{9.4} with $m=2$. Notice that as $p$ increases, $\Delta u$ tends to be piecewise constant. This is an indication the solution satisfies the Poisson equation with piecewise constant right hand side albeit with an extremely complicated solution pattern that clearly warrants further investigation.
  }
  \begin{center}
    \subfigure[{\label{fig:a1}
        The approximation to $u$, the solution of the $4$-Bilaplacian.
    }]{
      \includegraphics[scale=\figscale,width=0.47\figwidth]{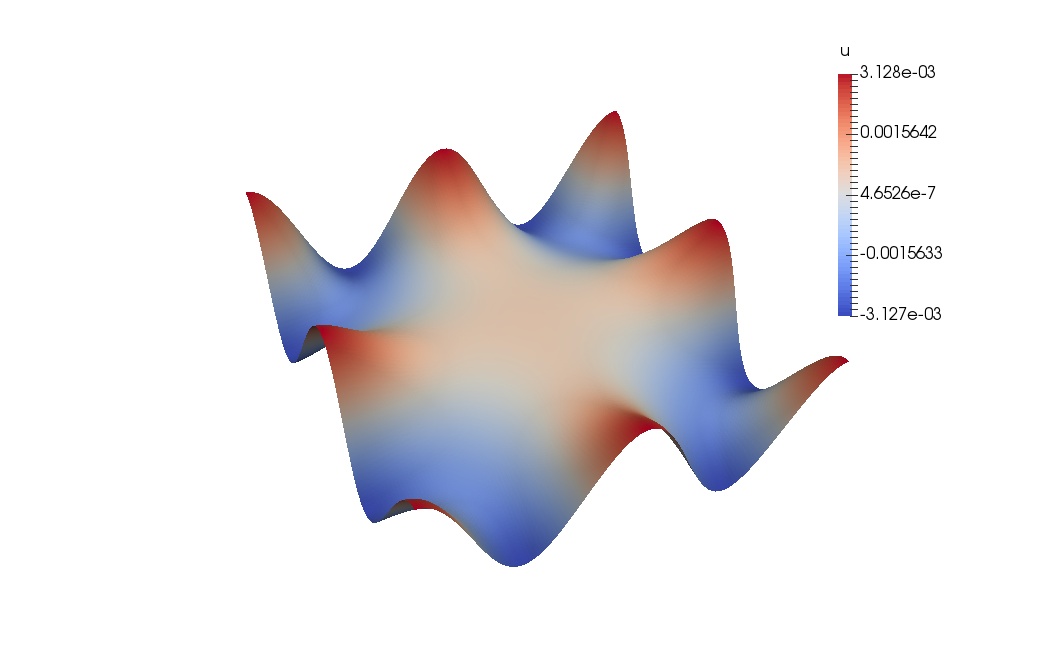}
    }  
    \hfill
    \subfigure[{\label{fig:a2}
        The approximation to $u$, the solution of the $142$-Bilaplacian.
    }]{
      \includegraphics[scale=\figscale,width=0.47\figwidth]{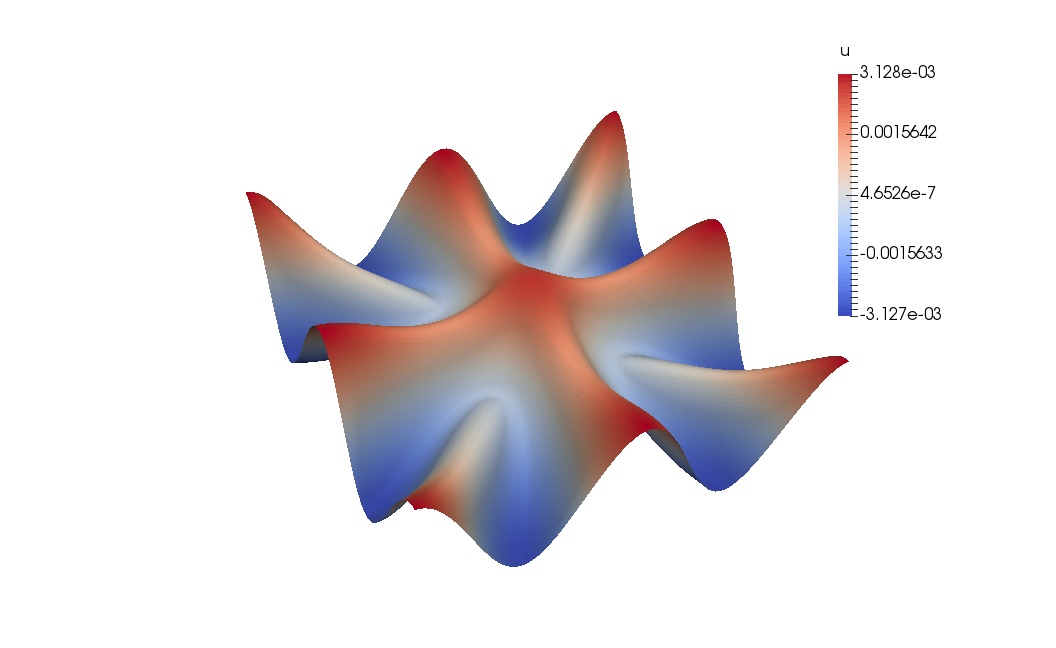}
    }
    \hfill
    \subfigure[{\label{fig:a1}
        The approximation to $\Delta u$, the Laplacian of the solution of the $4$-Bilaplacian.
    }]{
      \includegraphics[scale=\figscale,width=0.47\figwidth]{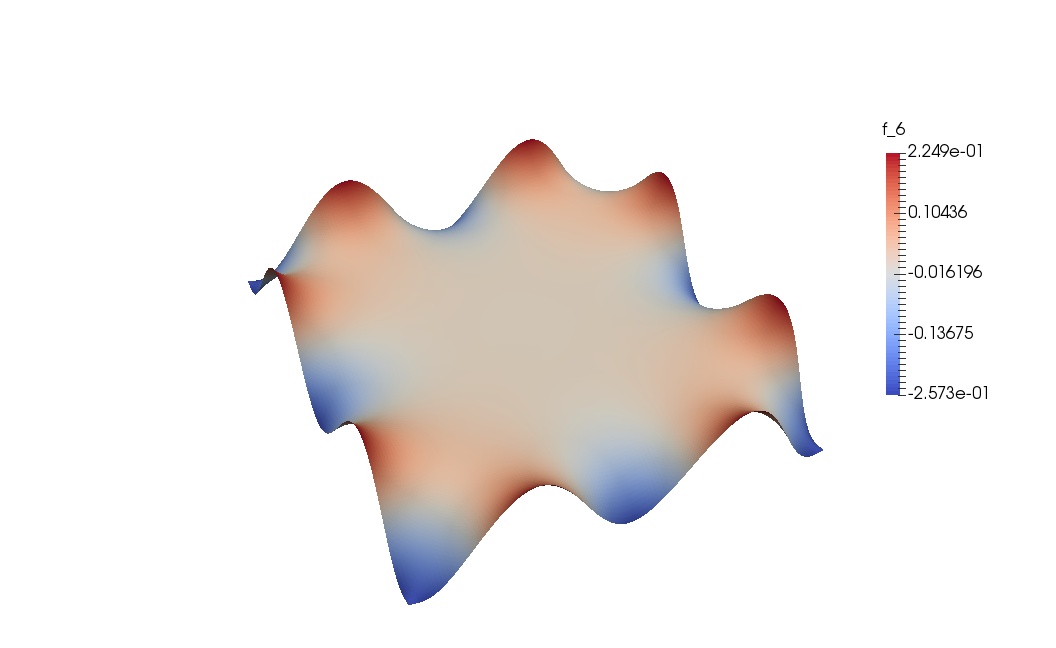}
    }
    \hfill
    \subfigure[{\label{fig:a2}
        The approximation to $\Delta u$, the Laplacian of the solution of the $42$-Bilaplacian.        
    }]{
      \includegraphics[scale=\figscale,width=0.47\figwidth]{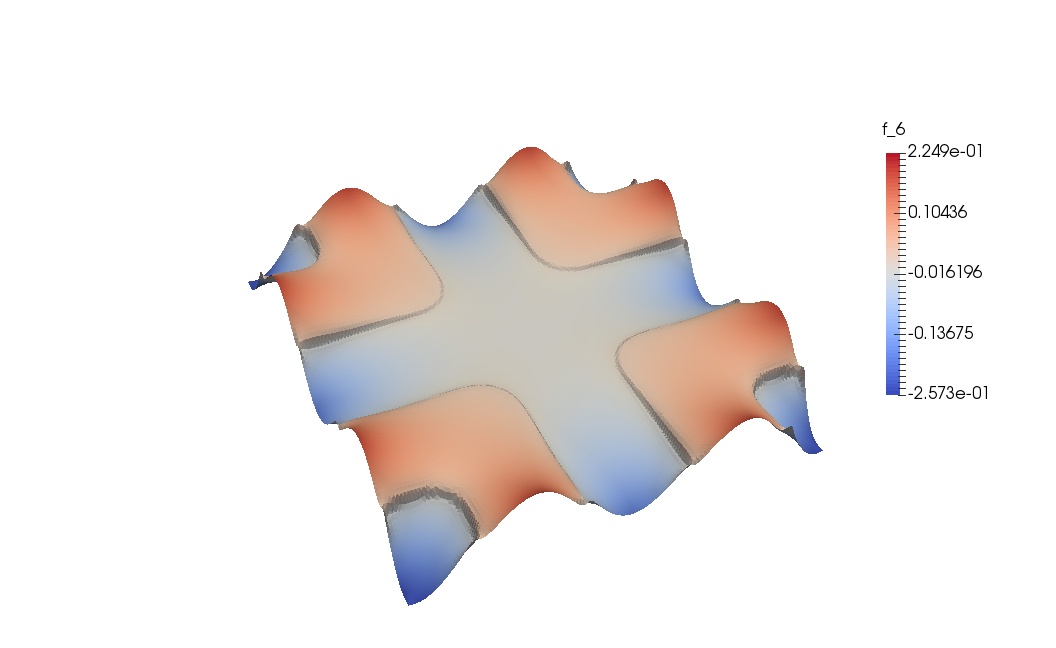}
    }
    \hfill
    \subfigure[{\label{fig:a2}
        The approximation to $\Delta u$, the Laplacian of the solution of the $68$-Bilaplacian.
    }]{
      \includegraphics[scale=\figscale,width=0.47\figwidth]{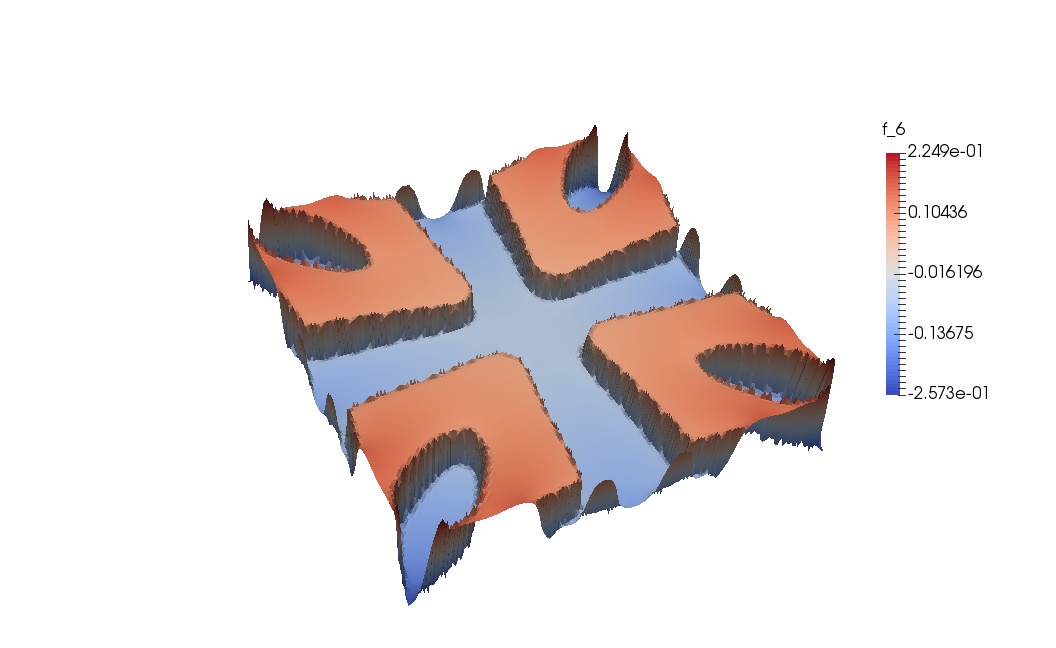}
    }
    \hfill
    \subfigure[{\label{fig:a2}
        The approximation to $\Delta u$, the Laplacian of the solution of the $142$-Bilaplacian.
    }]{
      \includegraphics[scale=\figscale,width=0.47\figwidth]{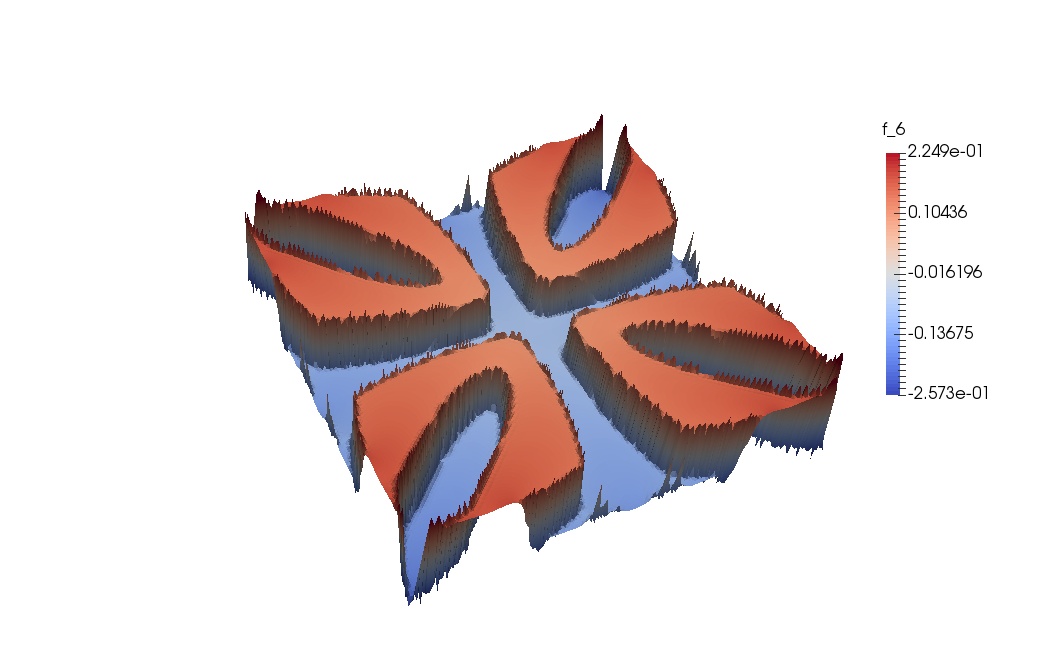}
    }
 
  \end{center}
\end{figure}

\begin{figure}[!ht]
  \caption[]
  {\label{fig:floppydonkeydick5}
    {\bf Test 3c:} A mixed finite element approximations to an $\infty$-Biharmonic function using $p$-Biharmonic functions for various $p$ for the problem given by \eqref{eq:2dpbiharm} and \eqref{9.4} with $m=3$. Notice that as $p$ increases, $\Delta u$ tends to be piecewise constant. This is an indication the solution satisfies the Poisson equation with piecewise constant right hand side albeit with an extremely complicated solution pattern that clearly warrants further investigation.
  }
  \begin{center}
    \subfigure[{\label{fig:a1}
        The approximation to $u$, the solution of the $4$-Bilaplacian.
    }]{
      \includegraphics[scale=\figscale,width=0.47\figwidth]{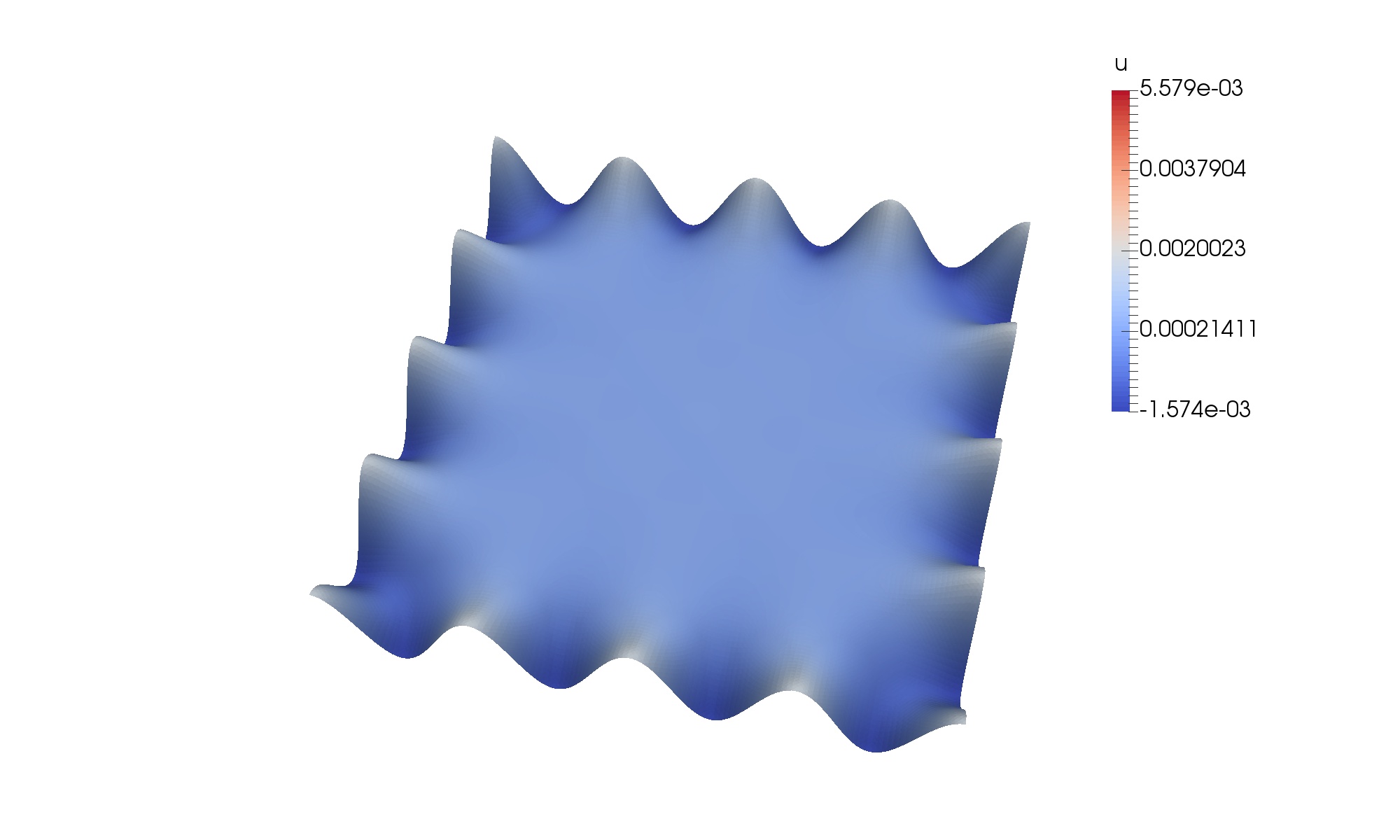}
    }  
    \hfill
    \subfigure[{\label{fig:a2}
        The approximation to $u$, the solution of the $142$-Bilaplacian.
    }]{
      \includegraphics[scale=\figscale,width=0.47\figwidth]{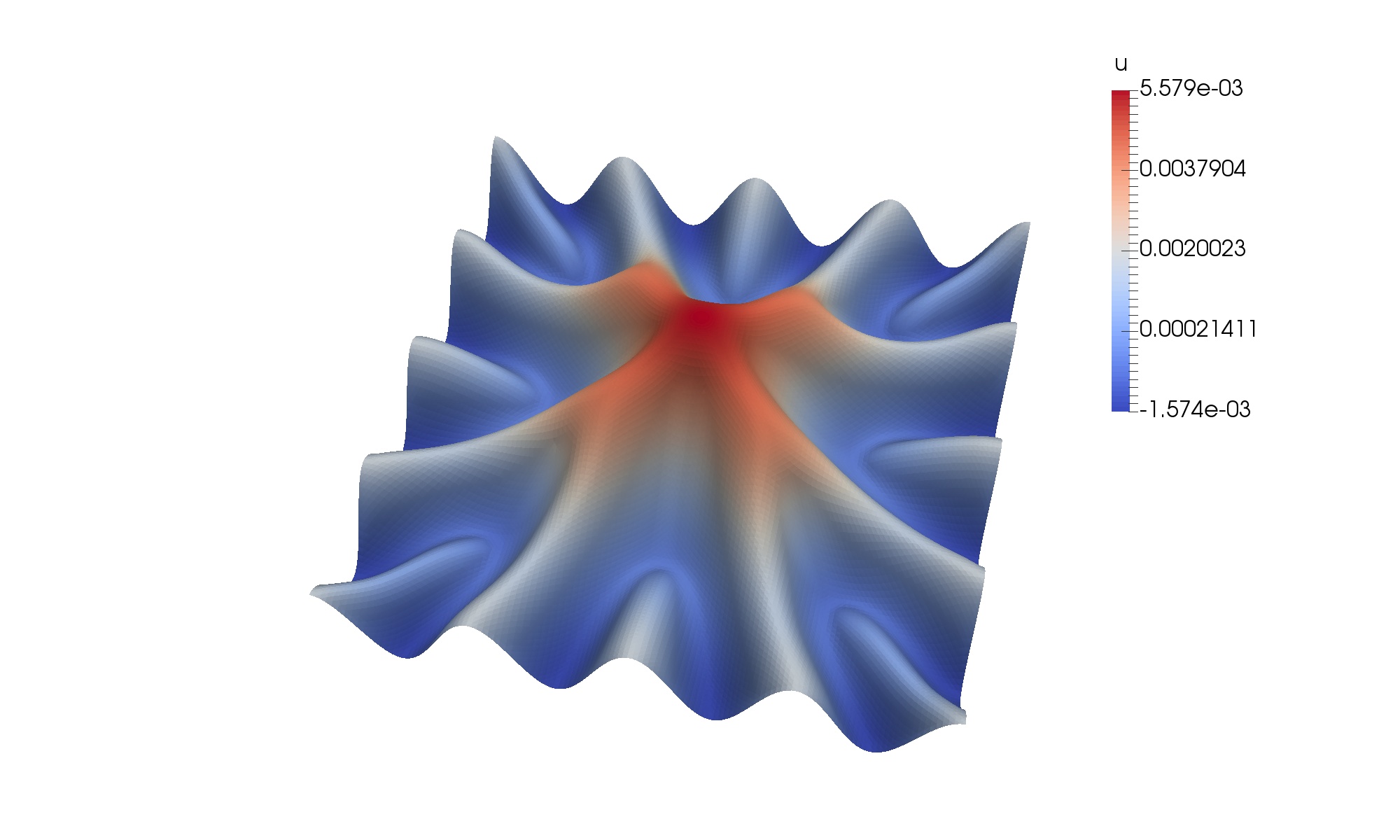}
    }
    \hfill
    \subfigure[{\label{fig:a1}
        The approximation to $\Delta u$, the Laplacian of the solution of the $4$-Bilaplacian.
    }]{
      \includegraphics[scale=\figscale,width=0.47\figwidth]{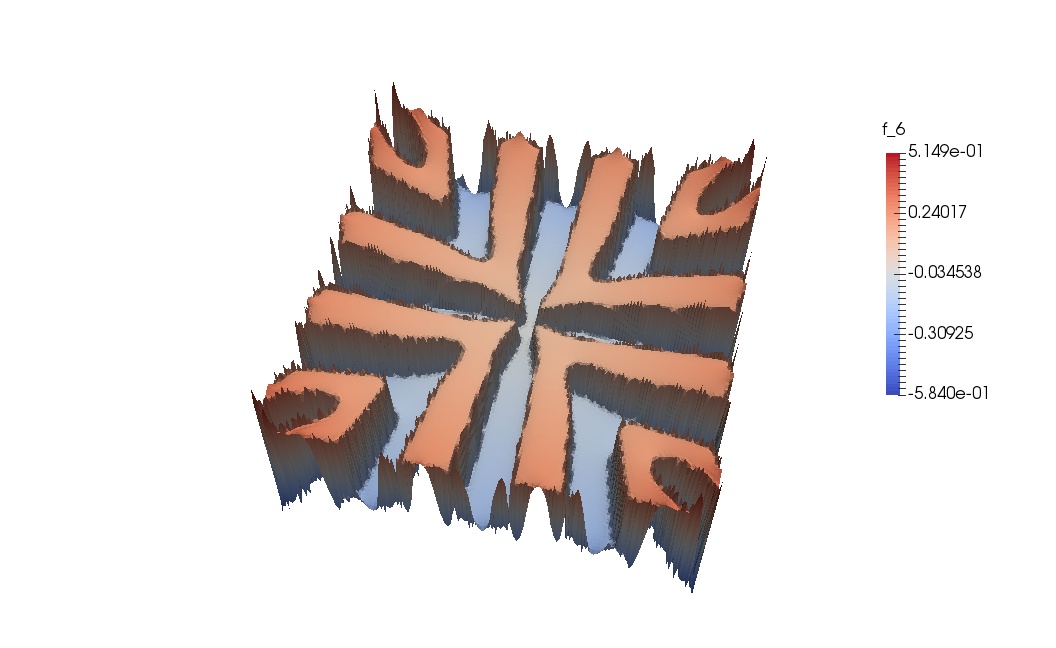}
    }
    \hfill
    \subfigure[{\label{fig:a2}
        The approximation to $\Delta u$, the Laplacian of the solution of the $42$-Bilaplacian.        
    }]{
      \includegraphics[scale=\figscale,width=0.47\figwidth]{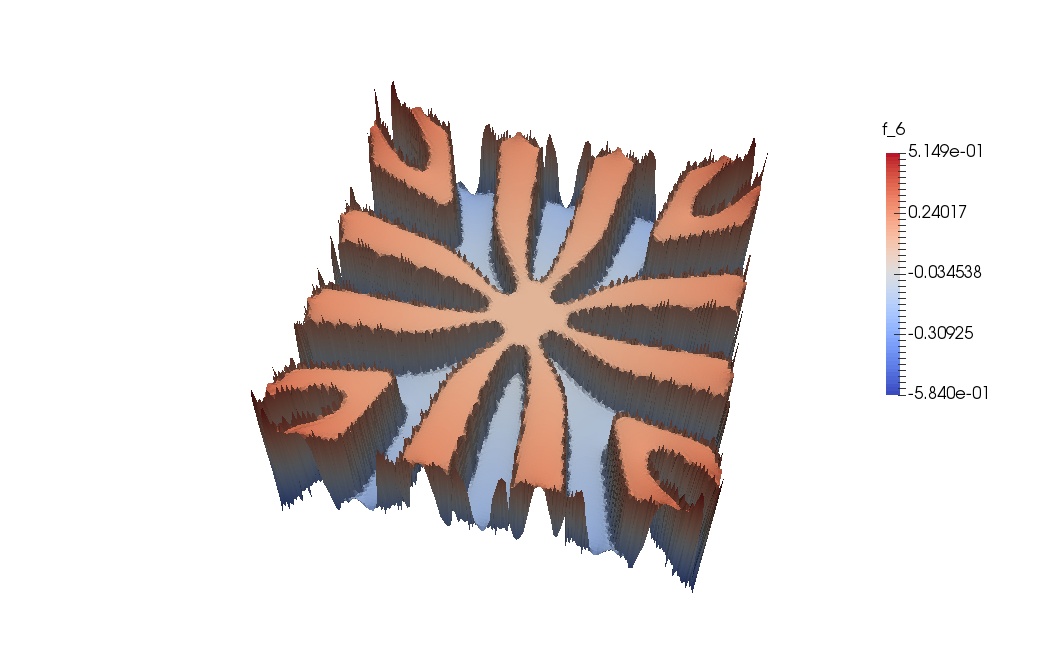}
    }
    \hfill
    \subfigure[{\label{fig:a2}
        The approximation to $\Delta u$, the Laplacian of the solution of the $68$-Bilaplacian.
    }]{
      \includegraphics[scale=\figscale,width=0.47\figwidth]{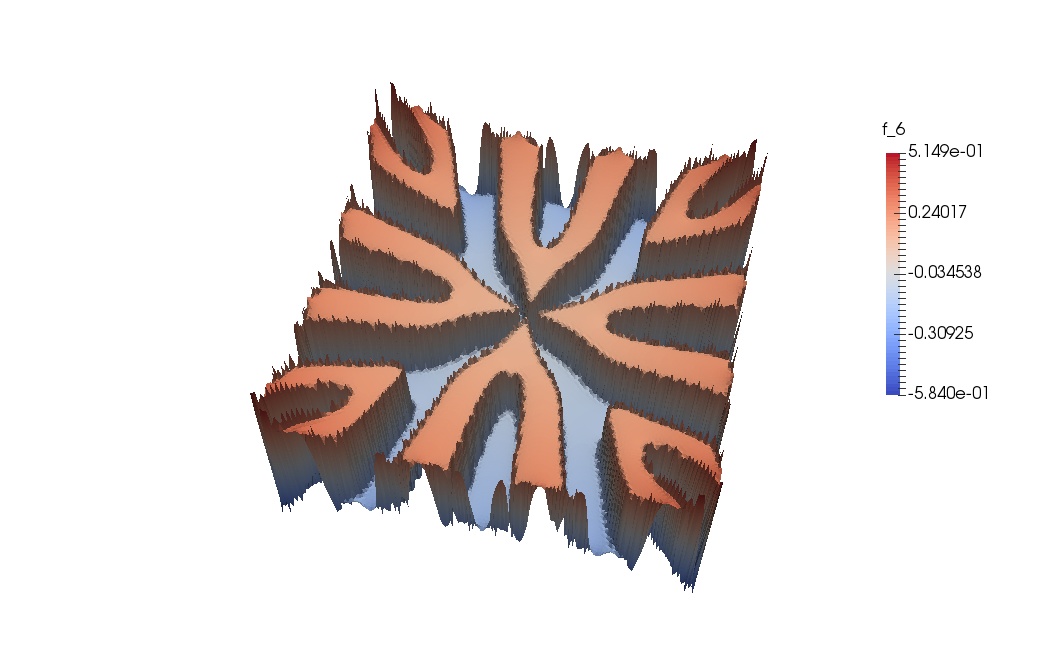}
    }
    \hfill
    \subfigure[{\label{fig:a2}
        The approximation to $\Delta u$, the Laplacian of the solution of the $142$-Bilaplacian.
    }]{
      \includegraphics[scale=\figscale,width=0.47\figwidth]{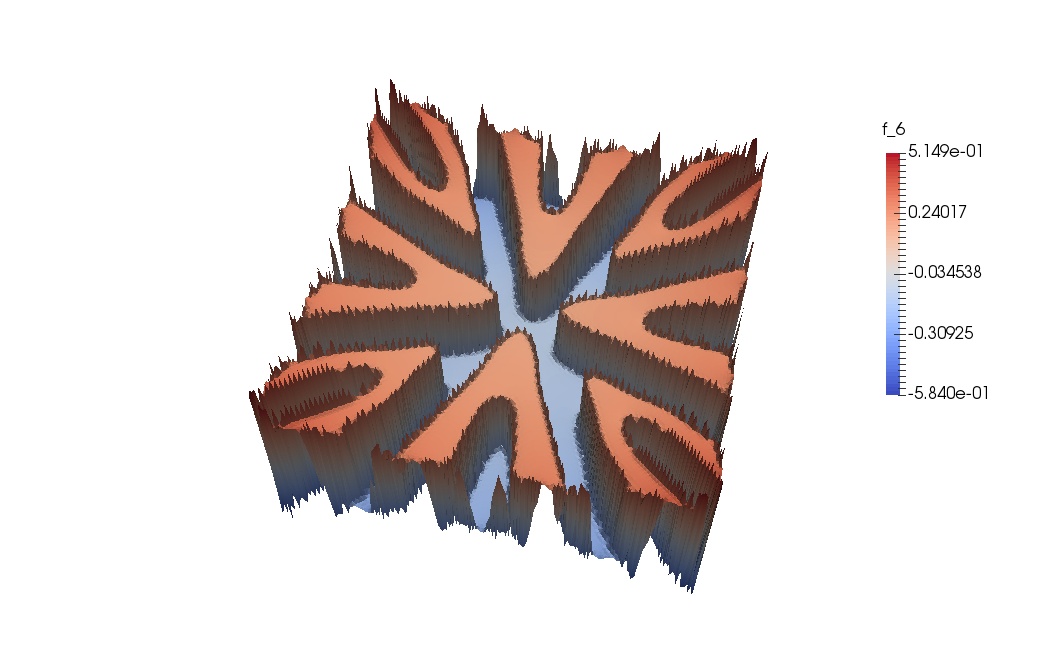}
    }
 
  \end{center}
\end{figure}

Now we illustrate some of the complicated behaviour of the
$p$-Bilaplacian for $d=2$:
\begin{equation}
  \label{eq:2dpbiharm}
  \left\{ \ \ \
  \begin{array}{rl}	
    \Delta \big( |\Delta u|^{p-2} \Delta u \big) \,=\,0, \ \ & \text{ in }\W = [-1,1]^2,
    \smallskip\\
    u \,= \,g, \ \ &   \text{ on }\partial\W,
    \smallskip \\
    \D u \,=\, \D g, & \text{ on }\partial\W,
  \end{array}
  \right.
\end{equation}
where $g$ is prescribed as
\begin{equation}
  \label{9.4}
  g(x,y) = \tfrac{1}{m 20} \cos{m \pi x} \cos{m \pi y},
\end{equation}
for various values of $m$.
We simulate the $p$-Bilaplacian for increasing values of $p$ and present the results in Figures \ref{fig:floppydonkeydick3},\ref{fig:floppydonkeydick4} and \ref{fig:floppydonkeydick5} indicating that in the limit the $\infty$-Biharmonic function should be piecewise quadratic however the behaviour is quite unexpected and complicated interface patterns emerge even with this relatively simple boundary data.

\section{Conclusion}

In this work we constructed a numerical method for the approximation
of solutions of the $p$-Bilaplacian equation. We were able to
analytically show convergence of the numerical approximation and, in
particular, to the solution of the limiting problem of the
$\infty$-Bilaplacian. This is particularly challenging as it is a
third order fully nonlinear PDE that is not in divergence form.

We have shown numerically that, for fixed $p$, our method converges
with rates that are better than the analysis predicted. This is well
documented in the case of similar lower order problems and can be
improved by using appropriate quasi-norms. We have utilised the
numerical method to make various interesting observations on the
structure of $\infty$-Biharmonic functions in that they are piecewise
quadratic over the domain with particularly complicated structures for
the interfaces.

{

\begin{figure}[!ht]
  \caption[]
  {\label{fig:floppydonkeydick3}
    {\bf Test 3a:} A mixed finite element approximations to an $\infty$-Biharmonic
    function using $p$-Biharmonic functions for various $p$ for the
    problem given by \eqref{eq:2dpbiharm} and \eqref{9.4} with
    $m=1$. Notice that as $p$ increases, $\Delta u$ tends to be
    piecewise constant. This is an indication the solution satisfies
    the Poisson equation with piecewise constant right hand side
    albeit with an extremely complicated solution pattern that clearly
    warrants further investigation.
  }
  \begin{center}
    \subfigure[{\label{fig:a1}
        The approximation to $u$, the solution of the $4$-Bilaplacian.
    }]{
      \includegraphics[scale=\figscale,width=0.47\figwidth]{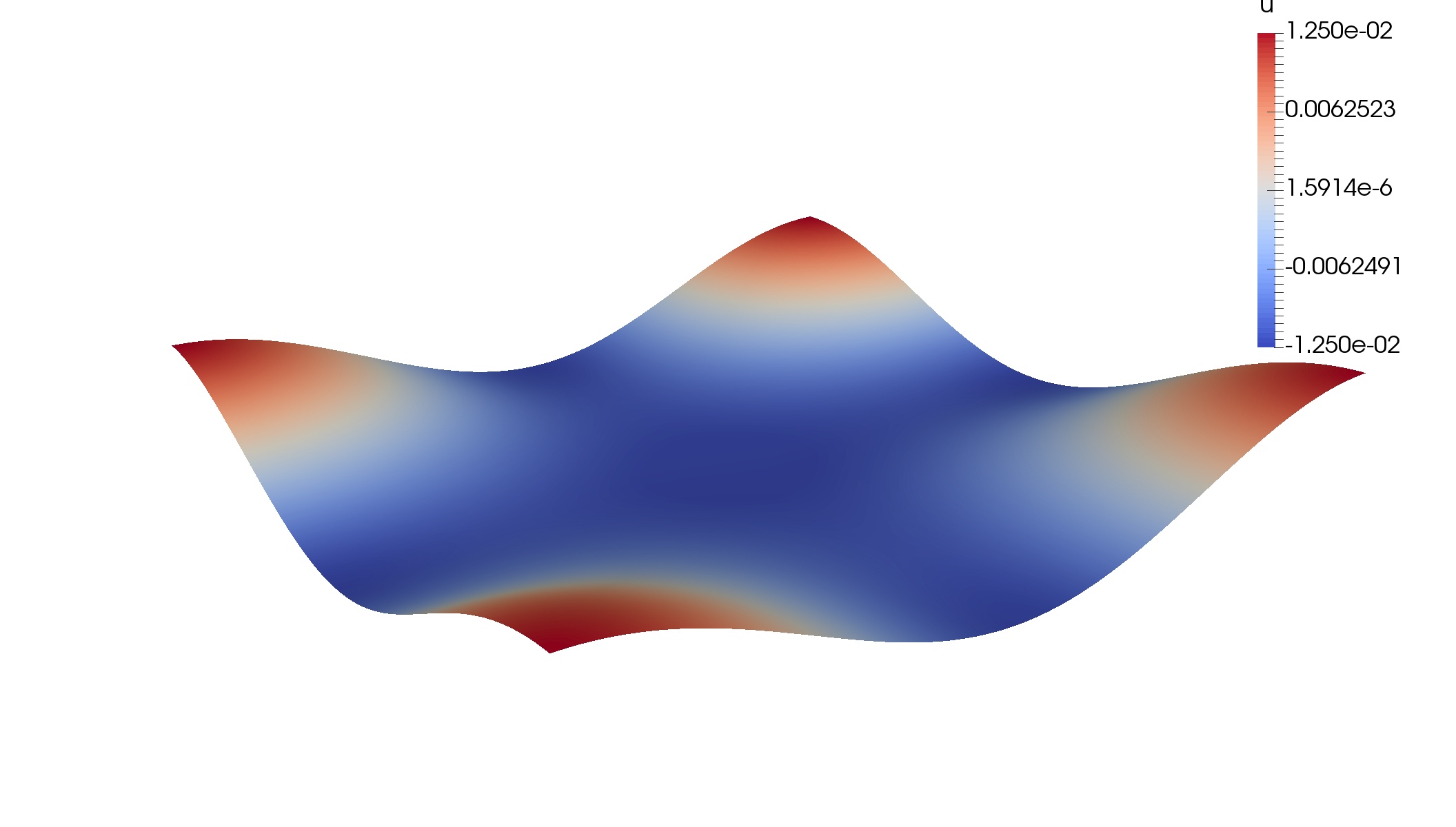}
    }  
    \hfill
    \subfigure[{\label{fig:a2}
        The approximation to $u$, the solution of the $142$-Bilaplacian.
    }]{
      \includegraphics[scale=\figscale,width=0.47\figwidth]{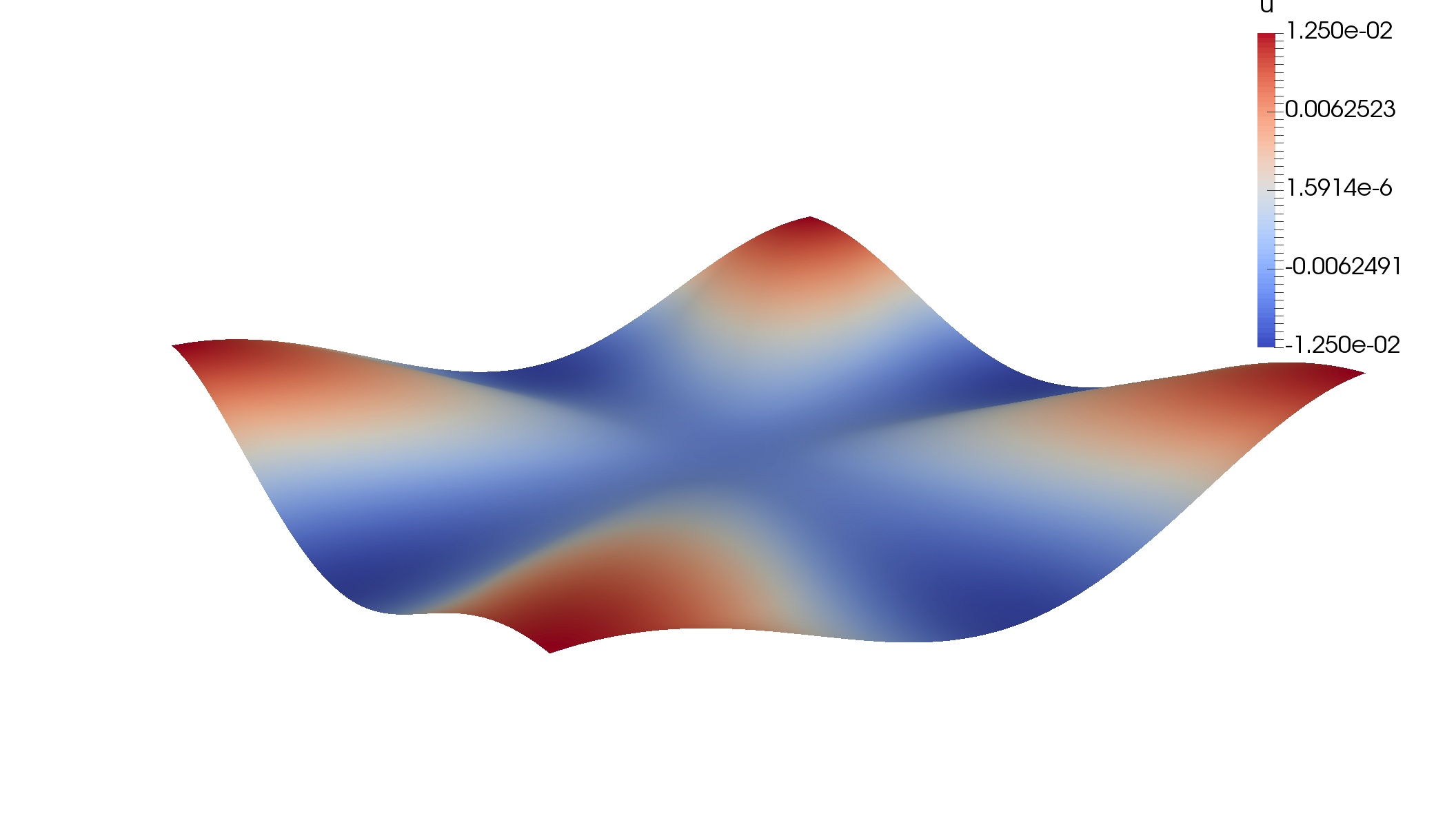}
    }
    \hfill
    \subfigure[{\label{fig:a1}
        The approximation to $\Delta u$, the Laplacian of the solution of the $4$-Bilaplacian.
    }]{
      \includegraphics[scale=\figscale,width=0.47\figwidth]{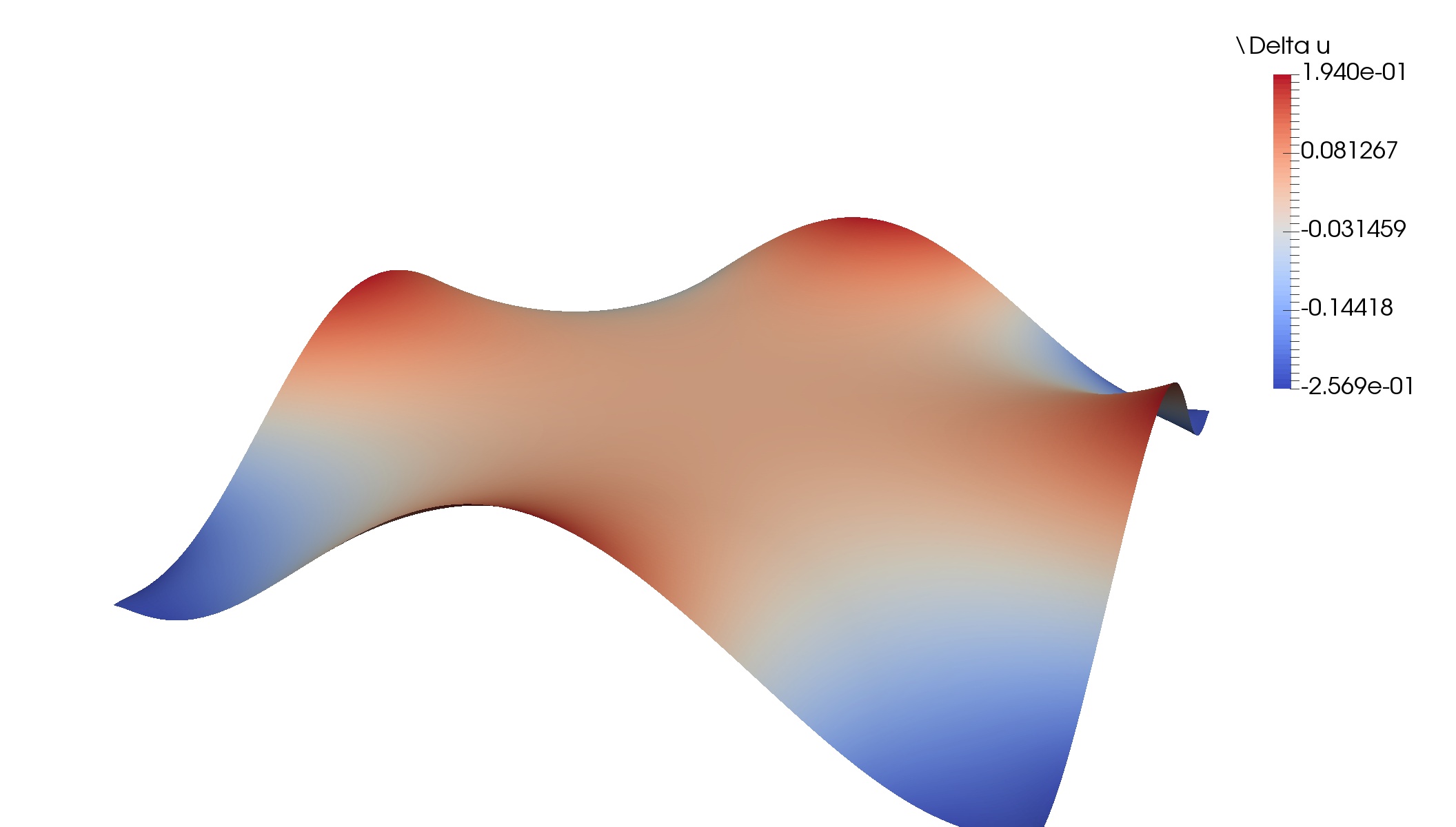}
    }
    \hfill
    \subfigure[{\label{fig:a2}
        The approximation to $\Delta u$, the Laplacian of the solution of the $42$-Bilaplacian.        
    }]{
      \includegraphics[scale=\figscale,width=0.47\figwidth]{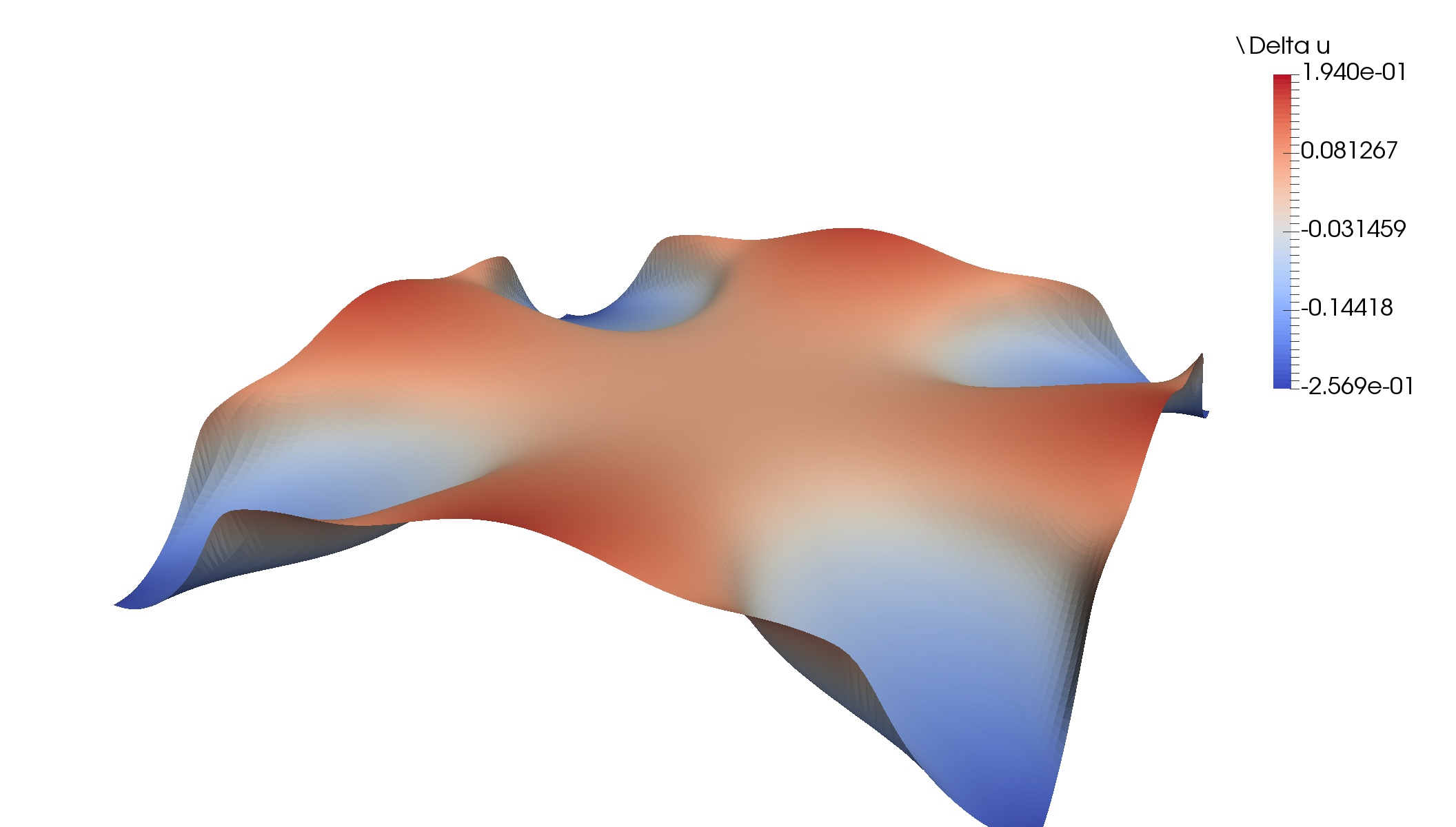}
    }
    \hfill
    \subfigure[{\label{fig:a2}
        The approximation to $\Delta u$, the Laplacian of the solution of the $68$-Bilaplacian.
    }]{
      \includegraphics[scale=\figscale,width=0.47\figwidth]{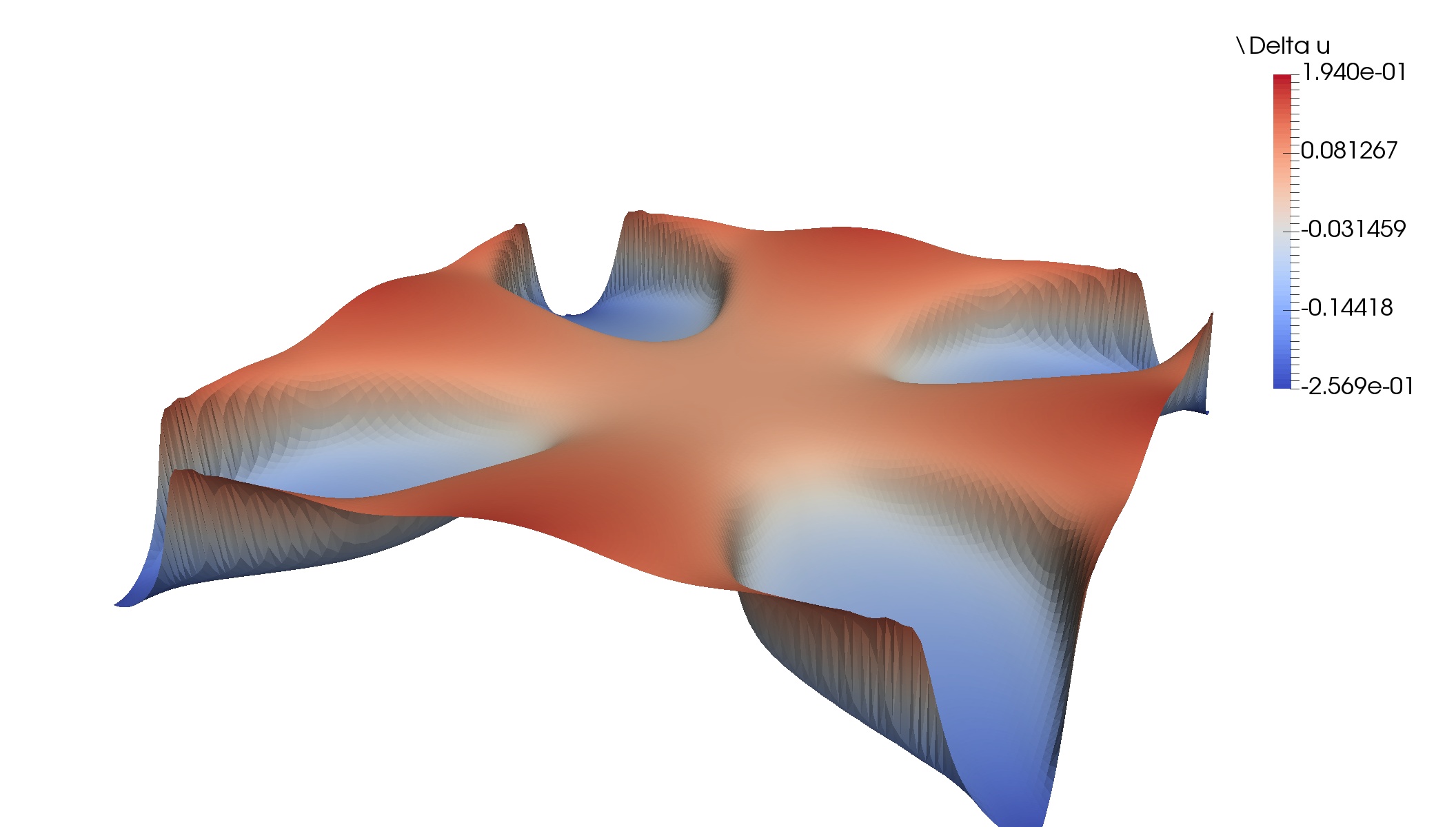}
    }
    \hfill
    \subfigure[{\label{fig:a2}
        The approximation to $\Delta u$, the Laplacian of the solution of the $142$-Bilaplacian.
    }]{
      \includegraphics[scale=\figscale,width=0.47\figwidth]{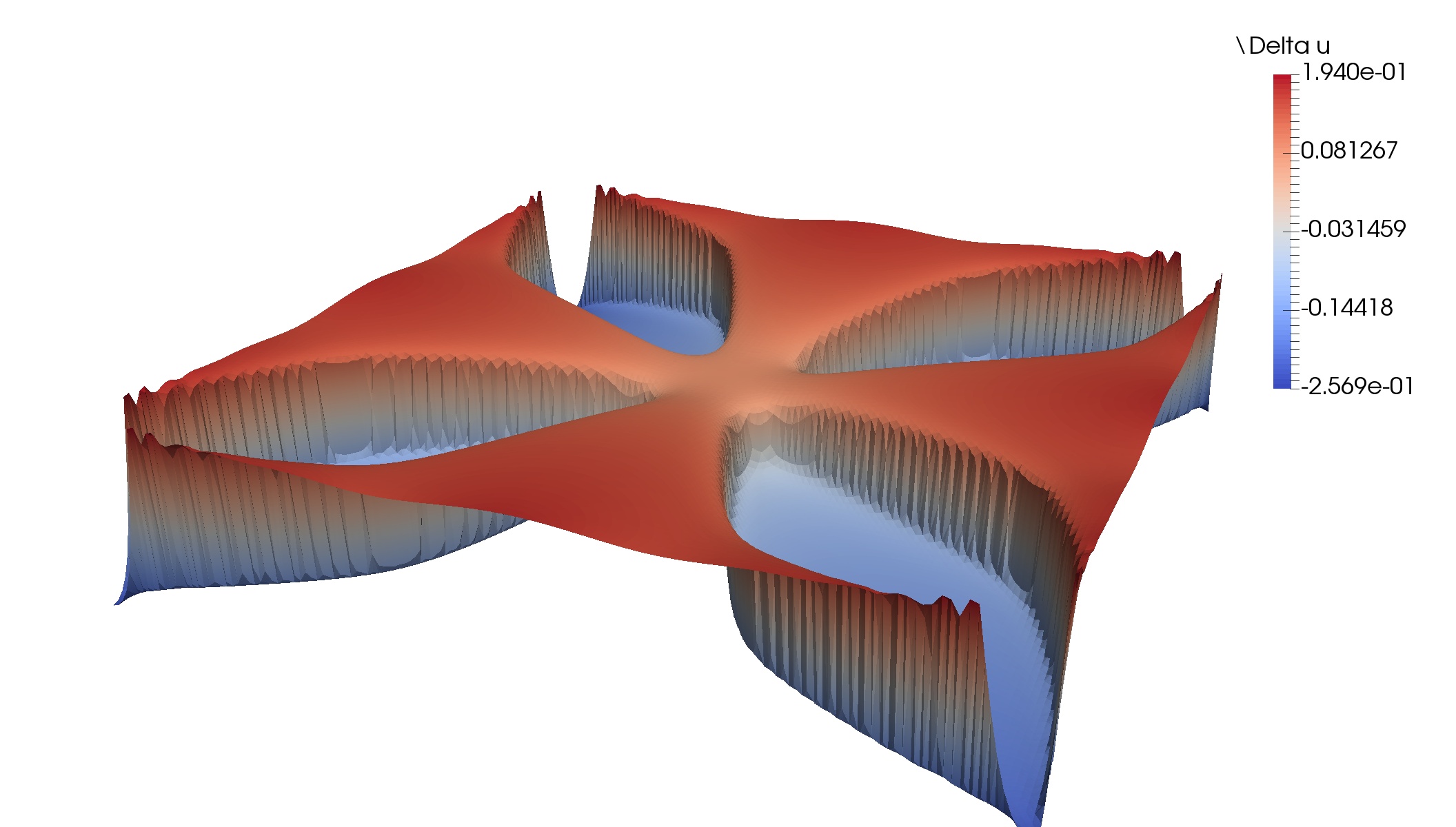}
    }
 
  \end{center}
\end{figure}

}

\printbibliography

\end{document}